\documentclass{article}
\usepackage[utf8]{inputenc}
\usepackage{titlesec}	
\usepackage[english]{babel}
\usepackage{pst-plot}
\usepackage{caption}
\usepackage{authblk}
\usepackage{subcaption}
\usepackage{amsthm,amssymb,amsmath,amsfonts,dsfont,calrsfs,upgreek,mathtools}
\usepackage{bbm}
\usepackage[mathscr]{euscript}
\usepackage{mathrsfs}

\newtheorem{mythm}{Theorem}[section]
\newtheorem{myexa}{Example}[section]
\newtheorem{mydef}{Definition}[section]

\newtheorem{mycor}{Corollary}[section]
\newtheorem{mypro}{Proposition}[section]

\newtheorem{myrem}{Remark}[section]

\newcommand{\A}{\mathcal{A}}
\newcommand{\B}{\mathcal{B}}
\renewcommand{\L}{\mathcal{L}}
\newcommand{\G}{\Gamma}
\newcommand{\tG}{\widetilde{\sigma}}
\newcommand{\id}{\mathds{1}}
\newcommand{\gid}{\mathbf{1}}
\newcommand{\idperm}{\mathrm{id}}
\newcommand{\pc}{\mathscr{P}_{\A}}
\newcommand{\g}{\textcircled{g}}
\newcommand{\Gal}{\mathbb{G}}
\newcommand{\Y}{\widetilde{Y}}

\usepackage{tikz}
\usetikzlibrary{automata, shapes, positioning, arrows, decorations.markings, decorations.pathreplacing, angles,quotes}
\tikzset{
    ->, % makes the edges directed
    >=stealth, % makes the arrow heads bold
    node distance=3cm, % specifies the minimum distance between two nodes. Change if necessary.
    every state/.style={thick, fill=gray!10}, % sets the properties for each ’state’ node
    every edge/.append style={line width=0.25mm}, % sets the properties for each ’state’ node
    initial text=$ $, % sets the text that appears on the start arrow
}
\title{Galois coverings of Schreier graphs of groups generated by bounded automata}
\author[1]{Hemant Bhate}
\author[2]{Daniele D'Angeli}
\author[1,3]{Asif Shaikh}
\author[4]{Dilip Sheth}

\affil[1]{Department of Mathematics, SP Pune University, Pune, India\\ \texttt{bhatehemant@gmail.com}}
\affil[2]{Dipartimento di Ingegneria, Università degli Studi Niccolò Cusano, Via Don Carlo Gnocchi, 3, Rome 00166, Italy\\ \texttt{daniele.dangeli@unicusano.it}}
\affil[3]{Department of Mathematics, R. A. Podar College of Commerce and Economics, Mumbai, India\\ \texttt{asif.shaikh@fulbrightmail.org}}
\affil[4]{Department of Mathematics, SP College, Pune, India\\ \texttt{shethdilip@yahoo.com}}

\date{\today}
\usepackage{graphicx}
\begin{document}

\maketitle

\begin{abstract}
We study coverings of Schreier graphs associated with groups generated by bounded automata. The action of such a group $G$ on the levels of the rooted tree $X^*$ produces a sequence of finite Schreier graphs $(\G_n)_{n \geq 0}$, each encoding the action of $G$ on words of length $n$. We show that the natural map $\G_{n+1} \to \G_n$ defines a $d$-sheeted covering, where $d = |X|$, and we characterize when this covering is Galois in terms of root permutations associated to the generating automaton. To analyze these coverings, we introduce a generalized replacement product of the graphs $\G_n$ and $\G_r$, which produces a graph isomorphic to $\G_{n+r}$. This construction generalizes the classical inflation procedure and provides a combinatorial model for the coverings that allows explicit computation of Ihara zeta and Artin $L$-functions in key examples. Our results highlight the interplay between automata structure, group actions, and graph coverings.
\end{abstract}


\section{Introduction}\label{sec:intro}

The notion of \emph{Schreier graphs} has become central to the study of group actions, symbolic dynamics, and the interface between algebra and geometry. Originating in the classical study of permutation representations of groups, these graphs arise from a triple $(G, H, S)$ consisting of a group $G$, a subgroup $H$, and a symmetric generating set $S$. They generalize Cayley graphs and provide a flexible and algebraically rich way to represent the action of a group on a coset space. The term \emph{Schreier graph} first appeared in the monograph by Coxeter and Moser~\cite{coxeter2013generators}, and has since become a key object of investigation in modern group theory and dynamics.

Remarkably, every regular graph of even degree can be realized as a Schreier graph, allowing one to apply group-theoretic and combinatorial techniques to their study. Schreier graphs have proved instrumental in addressing deep questions across mathematics. For instance, V.~Nekrashevych recently resolved R.~Grigorchuk's question on the existence of a finitely generated simple group of intermediate growth by analyzing the linear repetitiveness of associated Schreier graphs~\cite{nekrashevych2018palindromic}. Moreover, Schreier graphs exhibit striking spectral properties. The first examples of regular graphs with Cantor spectrum were obtained in this framework, as well as the first example of a group with discrete spectral measure, thereby answering a question of Atiyah (see~\cite{bartholdi2003fractal} and the references therein).

An especially rich class of Schreier graphs arises from groups generated by \emph{bounded automata}, which act level-transitively on rooted trees and yield self-similar actions. Famous examples include the Grigorchuk group, the Gupta--Sidki $p$-group, the Fabrykowski--Gupta group, and the Brunner--Sidki--Vieira torsion-free group (BSV). These groups exhibit unusual properties such as intermediate growth, amenability without subgroup stability, and atypical spectral behavior, many of which are reflected in the structure of their Schreier graphs.

In this paper, we study coverings of Schreier graphs arising from groups generated by bounded automata. Such groups act faithfully and level-preservingly on the regular rooted tree $X^*$, where $X$ is the input alphabet of size $d$. Each level $n$ corresponds to the set $X^n$ of words of length $n$, and the group acts on these levels via permutations. Given a symmetric generating set $S$, the Schreier graph $\G_n = \operatorname{Sch}(G, H_n, S)$ encodes the action of $G$ on $X^n$, where $H_n$ is the stabilizer in $G$ of a vertex of level $n$. This yields a sequence of finite Schreier graphs $(\G_n)_{n \geq 0}$ capturing the level-wise dynamics of $G$.

There are natural maps $\G_{n+1} \to \G_n$, and understanding when these maps are coverings, and further, when they are Galois coverings, forms a central theme of this work. In the case of a particular torsion-free weakly branch group defined by a three-state automaton, Grigorchuk and Żuk showed in Proposition~8.1 of~\cite{grigorchuk2002torsion} that $\G_{r+1}$ is a covering of $\G_r$. To study such maps systematically, we recall (see Section~\ref{sec:preliminaries}) that a graph covering is a locally bijective morphism between graphs.

We generalize this observation to the full class of groups generated by bounded automata. We prove that $\G_{n+r}$ is a $d^n$-sheeted unramified covering of $\G_r$. Moreover, we provide a characterization of when this covering is Galois. To analyze these coverings further, we introduce a new construction called the \emph{generalized replacement product}, which describes $\G_{n+r}$ explicitly as the product $\G_n \g \G_r$. This construction generalizes the classical inflation process introduced by Nekrashevych \cite{nekrashevych2005self}, refines the replacement product viewpoint, and provides a combinatorial model for these coverings. 

Using this framework, we provide a characterization of the conditions under which such coverings are Galois, meaning that there exist $d$ automorphisms of the covering graph that preserve the covering map and together form a group of order $d$. To this end, we construct automorphisms $\sigma$ of $\G_{n+1}$ arising from root permutations of reachable states. These automorphisms preserve adjacency and the covering structure, and under suitable conditions, generate the Galois group of the covering.

In the second part of the paper, we investigate the \emph{Ihara zeta} functions and \emph{Artin $L$}--functions associated with these Schreier graphs. These zeta functions encode refined information about the cycle structure of the graphs and its symmetry, and they provide a natural bridge between combinatorial, spectral, and representation-theoretic aspects of graph coverings. While Grigorchuk and {\.Z}uk ~\cite{grigorchuk2004ihara} studied the Ihara zeta functions for \emph{infinite} Schreier graphs of these groups, our focus is on the corresponding theory for the \emph{finite} graphs appearing at each level of the covering sequence. Exploiting the Galois structure of these coverings, we express the Ihara zeta functions in terms of Artin $L$--functions and carry out explicit computations in several examples. 

We also highlight connections to concrete combinatorial constructions. A notable example is the Tower of Hanoi group $\mathbb{T}^{(k)}$, generated by a finite automaton over an alphabet of size $k$, which models the classical puzzle on $k$ pegs. The associated Schreier graphs $\G_r^k$ form a sequence arising from the group action on the levels of a rooted tree. Interestingly, this group does not satisfy the conditions of our Galois covering criterion. Our results thus imply that the natural covering $\G_{r+1}^k \to \G_r^k$ is not normal. This example demonstrates both the sharpness and the applicability of our characterization. 

Before we close the introduction, we briefly outline the structure of the remainder of the paper. Section~\ref{sec:preliminaries} introduces the necessary background on graph coverings, automaton groups, and Schreier graphs. Section~\ref{sec:inflation} recalls the inflation process for Tile graphs. Section~\ref{sec:replacement} introduces the generalized replacement product and shows how it relates to the graph $\G_{n+r}$. Section~\ref{sec:galois} presents our main characterization theorem for Galois coverings in terms of automaton structure. In Section~\ref{sec:zeta}, we apply this to compute Ihara zeta and Artin $L$-functions for Schreier graphs of various automaton groups. Finally, Section~\ref{sec:nonbounded} considers examples beyond bounded automata and concludes with an open question.

\section{Preliminaries}\label{sec:preliminaries}

This section lays the foundational terminology and tools used throughout the paper. We begin with the notion of graph coverings, with particular emphasis on Galois (or normal) coverings, and then introduce groups generated by bounded automata together with their associated Schreier graphs. The terminology and conventions related to graph coverings and the corresponding zeta functions are largely adopted from \cite{terras2010zeta}, while our treatment of groups generated by bounded automata primarily follows \cite{bondarenko2007groups,bondarenko2016ends}.

\subsection{Galois coverings of graphs}

All graphs in this paper are finite, connected, labeled and may have multiple edges and loops. A \emph{graph morphism} $\pi: \tG \to \G$ is a map between graphs that preserves adjacency.

A morphism $\pi$ is a \emph{covering} if it is locally bijective: for every vertex $v \in \tG$, the map $\pi$ restricts to a bijection from the neighborhood of $v$ to the neighborhood of $\pi(v)$. We say that $\tG$ is an unramified covering of $\G$ if there exists a covering map $\pi: \tG \to \G$ that is surjective and satisfies the following condition: for every vertex $u \in \G$ and every $v \in \pi^{-1}(u)$, the neighborhood of $v$ maps bijectively onto the neighborhood of $u$. In particular, the size of the preimage (also called fiber) $\pi^{-1}(u)$ is constant for all $u \in \G$. If this number is $d$, we say that $\tG$ is a \emph{$d$-sheeted covering} of $\G$. Thus, the vertex set of $\tG$ can be represented as pairs $(u, i)$, where $u$ is a vertex of the graph $\G$ and $1 \leq i \leq d$.

\begin{mydef} \label{Def:Galois_cover}
A $d$-sheeted covering is called a \emph{normal} or \emph{Galois} covering if there exist $d$ automorphisms $\sigma: \tG \rightarrow \tG$ such that $\pi(\sigma(v)) = \pi(v)$ for all $v \in \G$. These automorphisms form a group called the Galois group, denoted $\Gal = Gal(\tG\mid\G)$.\end{mydef} 

We use $\gid$ to denote the identity element of the Galois group. When the Galois group is cyclic of order $d$, we say that $\tG$ is a $d$-fold cyclic covering of $\G$. After fixing a base sheet labeled $\gid$, the free and transitive action of $\Gal(\tG\mid\G)$ on each sheet allows us to canonically identify the $d$ sheets with the elements of the Galois group. We define the normalized Frobenius automorphism $\sigma(p) \in \Gal$ associated with a directed path $p$ of $\G$. This notion should be compared with the voltage assignment map of Gross and Tucker \cite{gross2001topological}.

\begin{mydef} \label{Def:nfauto}Suppose that $\tG\mid\G$ is normal with Galois group $\Gal = Gal(\tG\mid\G)$. For a path $p$ of $\G$, Proposition 13.3 of \cite{terras2010zeta} says there is a unique lift to a path $\tilde{p}$ of $\tG$ , starting on sheet $\gid$, having the same length as $p$. If $\tilde{p}$ has its terminal vertex on the sheet labeled $g \in \Gal$, define the normalized Frobenius automorphism $\sigma(p) \in \Gal$ by $$\sigma(p) \coloneqq g.$$
\end{mydef}
In this paper, we are primarily interested in such coverings arising from the action of bounded automaton groups on levels of regular rooted trees.

\subsection{Groups generated by bounded automata}
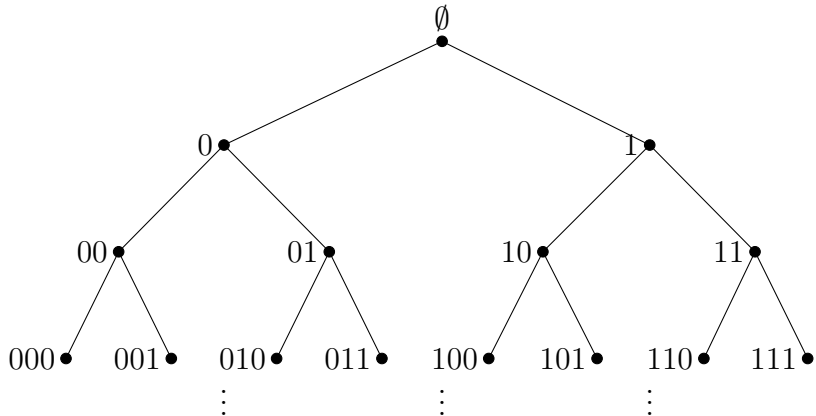
\begin{figure}[!htb]
	\centering
	\begin{tikzpicture}[scale=1.5]
	\draw[-] (2.50,2.10)  to  (1.05,1.41) to (0.35,0.70) to (0.00,-0.01);
	\draw[-] (2.50,2.10)  to  (3.88,1.41) to (4.58,0.70) to (4.93,-0.01);
	\draw[-] (1.05,1.41)  to  (1.75,0.70) to (2.10,-0.01);
	\draw[-] (0.35,0.70) to (0.70,-0.01);
	\draw[-] (1.75,0.70) to (1.40,-0.01);
	\draw[-] (3.88,1.41) to (3.17,0.70) to (2.81,-0.01);
	\draw[-] (3.17,0.70) to (3.53,-0.01);
	\draw[-] (4.58,0.70) to (4.24,-0.01);
	\draw[fill] (2.50,2.10) circle [radius=0.035] node[above ]  {$\emptyset$};
	\draw[fill] (1.05,1.41) circle [radius=0.035] node[left]  {$0$};
	\draw[fill] (3.88,1.41) circle [radius=0.035] node[left=5 pt]  {$1$};
	\draw[fill] (0.35,0.70) circle [radius=0.035] node[left]  {$00$};
	\draw[fill] (1.75,0.70) circle [radius=0.035] node[left]  {$01$};
	\draw[fill] (3.17,0.70) circle [radius=0.035] node[left]  {$10$};
	\draw[fill] (4.58,0.70) circle [radius=0.035] node[left]  {$11$};
	\draw[fill] (0.00,-0.01) circle [radius=0.035] node[left]  {$000$};
	\draw[fill] (0.70,-0.01) circle [radius=0.035] node[left]  {$001$};
	\draw[fill] (1.40,-0.01) circle [radius=0.035] node[left]  {$010$};
	\draw[fill] (2.10,-0.01) circle [radius=0.035] node[left]  {$011$};
	\draw[fill] (2.81,-0.01) circle [radius=0.035] node[left]  {$100$};
	\draw[fill] (3.53,-0.01) circle [radius=0.035] node[left]  {$101$};
	\draw[fill] (4.24,-0.01) circle [radius=0.035] node[left]  {$110$};
	\draw[fill] (4.93,-0.01) circle [radius=0.035] node[left]  {$111$};
    \node at (2.50,-0.25) {$\vdots$};
	\node at (1.05,-0.25) {$\vdots$};
	\node at (3.88,-0.25) {$\vdots$};
	\end{tikzpicture}
	\caption{Binary rooted tree.}
	\label{Figure:binarytree}
\end{figure}

Let $X = \{x_1,\dots, x_d\}$ be a finite alphabet with $d = |X| > 1$. The set $X^*$ of all finite words over $X$ (including the empty word $\emptyset$) forms a rooted $d$-ary tree with root $\emptyset$ and edges connecting $v$ to $vx$ for $v \in X^*$ and $x \in X$. The $n$-th level of the tree is the subset $X^n$. Figure~\ref{Figure:binarytree} illustrates the rooted tree $X^*$ in the case $d = 2$.

\begin{mydef}\label{def:self-similar-gr} Let $Aut(X^*)$ denote the group of automorphisms of this rooted tree. A group $G \leq Aut(X^*)$ is said to act self-similarly if for each $g \in G$ and $x \in X$, there exists $g|_x \in G$ such that
$$g(xv) = g(x)\, g|_x(v), \qquad \text{for all } v \in X^*.$$  
\end{mydef}
The element $g|_x$ is called the \emph{restriction} of $g$ at $x$. Inductively,
$$g|_{x_n \cdots x_2 x_1}= g|_{x_n}|_{x_{n-1}} \cdots |_{x_1}, \quad x_n \cdots x_2 x_1 \in X^*.$$
These satisfy the identities
\begin{equation}
g(uv) = g(u)g|_u(v), \quad g|_{uv} = g|_u|_v, \quad (gh)|_u = g|_{h(u)}\cdot h|_u. \label{eqn:restriction_g}
\end{equation}
Self-similarity gives rise to a recursive description: $G \leq Sym(X)\wr G$, where $\wr$ denotes the wreath product. The associated \emph{wreath recursion} map
$$\alpha: G \to Sym(X)\wr G, \qquad \alpha(g) = \psi_g \,(g|_{x_1},\dots,g|_{x_d})$$
records the \emph{root permutation} $\psi_g \in Sym(X)$ together with the restrictions $g|_{x_i}$.

\begin{mydef}\label{def:aut_gr}
  Let $\A=(S,X,t,o)$ be a finite automaton with state set $S$, input alphabet $X$, transition function $t:S\times X\to S$, and output function $o:S\times X\to X$. It is called \emph{invertible} if, for each $s\in S$, the map $o(s,\cdot):X\to X$ is a permutation. Every state $s\in S$ defines an automorphism of $X^*$ by the rule
$$s(xv) = o(s,x)\, (t(s,x)(v)), \qquad x\in X,\, v\in X^*.$$ Here the restriction $s|_x$ is exactly the next state $t(s,x)$, and the root permutation $\psi_s$ is $o(s,\cdot)$. Thus each state acts as a tree automorphism in a self-similar way, and the group $\langle \A \rangle$ generated by the states of $\A$ is a self-similar subgroup of $Aut(X^*)$, called an automaton group.
\end{mydef}
Conversely, a finitely generated self-similar group can be represented by an invertible finite automaton by choosing a generating set closed under restrictions. This viewpoint allows us to use the \emph{wreath recursion} to track root permutations of elements of $G$, while at the same time exploiting the \emph{automaton structure} to define post-critical sequences via the Moore diagram of $\A$.

Among automaton groups, those generated by bounded automata form a particularly important subclass, introduced by Sidki~\cite{sidki2000automorphisms}. include many fundamental examples and will be the primary focus of this work.

\begin{mydef}
    Let $g \in Aut(X^*)$. Define $\theta_k(g)$ as the number of words $v \in X^k$ such that $g|_v$ acts non-trivially on $X$. We say $g$ is \emph{bounded} if the sequence $\theta_k(g)$ is bounded in $k$. An invertible automaton is \emph{bounded} if all its states define bounded automorphisms.
\end{mydef}

This class of groups include many important examples like the Grigorchuk group, Gupta-Sidki group, and the Hanoi Tower groups, and they serve as the primary focus of this work. We use $\idperm$ to denote the identity permutation in the symmetric group $Sym(X)$, and $\id$ to denote the trivial state of the automaton.

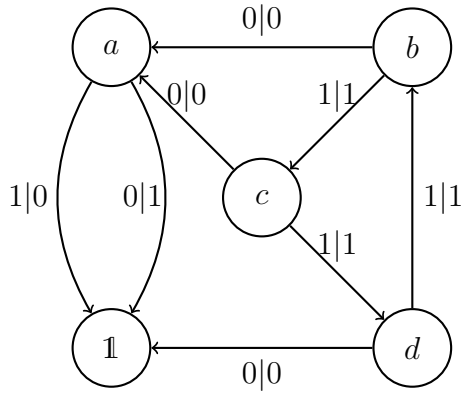
\begin{figure}[!htb]
	\centering
	\begin{tikzpicture}
	\tikzstyle{every state}=[draw=black, thick,fill=none]
    \node[state] at (-2,2) (a) {$a$};
	\node[state] at (2,2)  (b) {$b$};
	\node[state] at (2,-2) (d) {$d$};
	\node[state] at (-2,-2) (id) {$\id$};
	\node[state] at (0,0) (c) {$c$};
	\draw[thick,->] (b)  to  (a) ;
    \node[above] at (0,2) {$0|0$};
	\draw[thick,->] (b)  to  (c) ;
    \node[above] at (1,1) {$1|1$};
	\draw[thick,->] (d)  to  (b) ;
	\node[right] at (2,0) {$1|1$};
	\draw[thick,->] (d)  to  (id) ;
	\node[below] at (0,-2) {$0|0$};
	\draw[thick,->] (c)  to  (d) ;
	\node[above] at (1,-1) {$1|1$};
	\draw[thick,->] (c)  to  (a) ;
	\node[above] at (-1,1) {$0|0$};
	\draw[thick,->] (a)  to [out= 240, in=120] (id) ;
	\node[left=7pt] at (-1,0) {$0|1$};
	\draw[thick,->] (a)  to [out= 300, in=60] (id) ;
	\node[left] at (-2.7,0) {$1|0$};
	\end{tikzpicture}
	\caption{The automaton generating the Grigorchuk group.}
	\label{Figure:auto_grig}
\end{figure}

 \begin{myexa}Grigorchuk group \label{grig-group-ex} \\
	Let $X = \{0,1\}$. Define the automorphisms of the binary tree $X^*$ by the following wreath recursions:
\begin{equation}
\label{Equation:wr-grig-group}
a = \psi_a(\id,\id),\quad b = \psi_b (a, c),\quad c = \psi_c (a, d),\quad d = \psi_d (\id, b),
\end{equation}
where $\psi_a$ is the transposition $(0~1) \in Sym(X)$ and $\psi_b = \psi_c = \psi_d = \idperm$. 
 \end{myexa}   
    The group defined by the wreath recursions in \eqref{Equation:wr-grig-group} was constructed by R.~Grigorchuk as the first example of an infinite, periodic, finitely generated group. Moreover, it was the first known group with intermediate growth, meaning its growth function grows faster than any polynomial but slower than any exponential function \cite{grigorchuk1980burnside}. Since all sections in \eqref{Equation:wr-grig-group} are elements of the generating set, the group's action can be modeled by a finite invertible automaton. See Figure~\ref{Figure:auto_grig} for its Moore diagram.
 \begin{myexa} BSV torsion-free group \label{BSV-group-ex}\\
	Let $X = \{0, 1, \dots, d-1\}$. Define the automorphisms of the $d$-ary rooted tree $X^*$ by the following wreath recursions:
	\begin{equation} \label{Equation:wr-BSV-group}
	a = \psi_a(\id, \id, \cdots, \id, a),~ b = \psi_b(\id, \id, \cdots, \id, b^{-1}),~ b^{-1} = \psi_{b^{-1}}(\id, \id, \cdots, \id, b),
	\end{equation} where $\psi_a = (0~1~\cdots,d-1) = \psi_b, \psi_{b^{-1}} = \psi_b^{-1} \in Sym(X).$ 
 \end{myexa} 
This group was introduced by S. Sidki and E. F. Silva. (See \cite{sidki2001family}.) Note that {\it BSV torsion-free} $k$-group has transitive action on $d$-rooted tree, where $ 1 < d \in \mathbb{N}$. It is not difficult to show that {\it BSV torsion free group} generated by {\it bounded automaton}. This group is a generalization of a just-nonsolvable torsion-free group defined on the binary rooted tree by A. M. Bunner, S. Sidki and A. C. Vieira \cite{brunner1999just}.
\begin{myexa} Gupta-Sidki group.\\ Let $p =3$, so we have $X = \{0,1,2\}$. Define the automorphisms of the binary tree $X^*$ by the wreath recursions:
	\begin{equation}
	\label{Equation:wr-GS}
	a = \psi_a(b, b^{-1}, a),~ b = \psi_b(\id, \id, \id)
	\end{equation} 
	where $\psi_b = (0,1,2) \in Sym(X)$ and $\psi_a = \idperm$.
\end{myexa} The group is generated by $a, b, a^{-1}, b^{-1}$  and it was introduced by N. Gupta and S. Sidki. This group is an counter example to the general Burnside problem, which shows that a finitely generated group, all of whose elements have finite $p$-power order (for a fixed prime $p$), can be infinite. See \cite{gupta1983burnside}. 
\begin{myexa} Fabrykowski-Gupta group.\\ Let $X = \{0,1,2\}$. Define the automorphisms of the binary tree $X^*$ by the wreath recursions:
	\begin{eqnarray}
	\label{Equation:wr-FG}
	a = \psi_a(\id, \id, \id),~~ b = \psi_b(a,b, \id)
	\end{eqnarray}
	where $\psi_a = (0,1,2) \in Sym(X)$ and $\psi_b = \idperm$.    
\end{myexa} The group is generated by
$a, b$ and it was introduced by N. Gupta and J. Fabrykowski. See \cite{fabrykowski1987groups}. 
\begin{myexa} Basilica group \label{basilica-gr-ex}\\
	Let $X = \{0,1\}$. Define the automorphisms of the binary tree $X^*$ by the wreath recursions:
	\begin{equation} \label{Equation:wr-B-group}
	a = \psi_a (b, \id),~ b = \psi_b (a, \id),
	\end{equation}
	where $\psi_a = \idperm$ and $\psi_b = (0,~1) \in Sym(X).$    
\end{myexa} For generating automaton of Basilica group see Figure \ref{Figure:basilica-aut:left}. This group was introduced by R. Grigorchuk and A. {\.Z}uk \cite{grigorchuk2002torsion}. The Basilica group is related to the fractal called Basilica which is the Julia set of the polynomial $z^2-1$. 
 
The next section is dedicated to Schreier graphs and Tile graphs of groups generated by bounded automata.

\subsection{Schreier graphs of group actions}

Let $G$ be a group generated by a finite symmetric set $S$, and let $H$ be a subgroup of $G$. The \emph{Schreier graph} $Sch(G, H, S)$ is the directed labeled graph with vertex set $G/H$ (the set of right cosets), and edge set $G/H \times S$. For each coset $Hg$ and generator $s \in S$, there is a directed edge from $Hg$ to $Hgs$ labeled by $s$.

Alternatively, if $G$ acts on a set $M$ (on the right), the Schreier graph $Sch(G, M, S)$ has vertex set $M$ and a directed edge from $x$ to $xg$ labeled by $g$, for every $x \in M$ and $g \in S$. In this case, if the action is transitive and $m \in M$, then $Sch(G, M, S) \cong Sch(G, St_G(m), S)$, where $St_G(m)$ is the stabilizer of $m$ in $G$.

Let $G \leq Aut(X^*)$ be a group generated by a bounded automaton $\A$. The levels $X^n \subset X^*$ of the regular rooted tree are invariant under the action of $G$, and the group acts on each level $X^n$ via finite permutations. We define the $n$-level Schreier graph of this action as:
$$\G_n := Sch(G, X^n, S),$$
where $S$ is a finite symmetric generating set of $G$. The vertex set of $\G_n$ is $X^n$, and for each $x \in X^n$ and $s \in S$, there is a directed edge from $x$ to $s(x)$ labeled by $s$. This gives rise to a sequence of Schreier graphs $(\G_n)_{n \geq 0}$, connected via natural projection maps $\tau_n : X^{n+1} \to X^n$, which delete the last letter of a word. These induce graph morphisms $\tau_n : \G_{n+1} \to \G_n$, and thus the Schreier graphs form an inverse system.

We also consider a related subgraph called the \emph{Tile graph}, denoted by $\G_n'$, defined as follows.

\begin{mydef}\label{def:tile_gr}
The Tile graph $\G_n'$ has vertex set $X^n$, and an edge $\{v, s(v)\}$ is included if $s \in S$ and the section $s|_v = \id$.
\end{mydef}

Tile graphs are closely related to the combinatorics of fractals and appear naturally in the context of bounded automaton groups. They omit edges where the action depends non-trivially on deeper levels of the tree.

\begin{mydef} \label{def:post-crit-seq}  
A left-infinite sequence $\cdots x_2x_1$ over $X$ is called \emph{post-critical} if there exists a left-infinite path $\cdots e_2, e_1$ in the Moore diagram of $\A$, ending at a non-trivial state $s \in S$ and labeled by $\cdots x_2x_1|\cdots y_2y_1$ for some $y_i \in X$.
\end{mydef}

\begin{mydef} \label{def:reachable}
The non-trivial state $s$ in Definition \eqref{def:post-crit-seq} is referred to as a \emph{reachable state}.
\end{mydef}
Observe that, the existence of a post-critical sequence immediately implies the existence of a reachable state. Recall from Proposition 2.7 of \cite{bondarenko2016ends} that if $G$ is generated by a bounded automaton $\A$, then the set $\pc$ of post-critical sequences is finite. For a post-critical sequence $p = \cdots x_2x_1 \in \pc$, define $p_n = x_n \cdots x_2x_1$ as the associated post-critical vertex in $\G_n$ and $\G_n'$. If an edge $\{v, s(v)\}$ is present in $\G_n$ but missing in $\G_n'$, then $v$ is a post-critical vertex. The set of Tile graph edges is given by
\begin{align}
E(\G_n') &= E(\G_n) \setminus E, \label{eq:nonSchreier} \\
E &= \biggl\{ \{u, s(u)\} : u \in X^n,\; s \in S,\; s|_u \neq \id \biggr\}. \label{eqn:non-tile}
\end{align}

There is a natural one-to-one correspondence between the set $\pc$ and the set $P$ of post-critical vertices in the graphs $\G_n$ and $\G_n'$.
Before we conclude the preliminaries, let us state the assumptions that will remain in effect for the rest of the paper. Assumptions:
  \begin{enumerate}
    \item \label{ass:bdd}  $G$ is a group generated by a bounded automaton $\A$. Hence, the the set $\pc$ of post-critical sequences is non-empty and finite.
  \item The Tile graphs $\G_n'$ are connected for all $n$.
  \end{enumerate}
\section{Inflation of graphs}\label{sec:inflation}
In this section, we describe the process of \emph{inflation} in the finite graph-theoretical sense, following Bondarenko \cite{bondarenko2007groups}. This construction provides a combinatorial analog of the Tile diagram construction introduced by V. Nekrashevych in Section 3.9 of \cite{nekrashevych2005self}. It constructs higher-level Tile graphs $\G_{n+1}'$ from $\G_n'$ in a systematic manner. This construction relies on an auxiliary object called the model graph.

\begin{mydef} \label{def:model}
	The model graph $M = M_1$ has vertex set $V_M = \pc \times X$, and the edge set $E_M$ is defined as follows: 
\begin{eqnarray}
\label{Equation:model-edges}
E_M = \left\lbrace \begin{array}{l | l}
& \textnormal{ there exists a path } \cdots e_2e_1 \textnormal{ in the Moore }\nonumber\\ 
\{(p,x),(q,y)\} & \textnormal{ diagram of } \A, \textnormal{ which ends in the trivial  state }\\
&  \textnormal{ and which is labelled by where }  px|qy, p,q \in \mathcal{P_A}\\
& \textnormal{ and } ~x,y \in X \end{array} \right\rbrace.
\end{eqnarray}
\end{mydef}
\begin{myrem}
A vertex $(p, x) \in \pc \times X$ in the model graph $M$ is called a post-critical vertex if the corresponding sequence $px$ is post-critical. Note that there is a one-to-one correspondence between post-critical vertices of $M$ and elements of $\pc$.
\end{myrem}

We now recall Theorem 2 from \cite{bondarenko2016ends}, which formally describes how to obtain the next-level Tile graph $\G_{n+1}'$ via inflation from $\G_n'$ using the model graph.

\begin{mythm}
\label{Theorem:bon-inflation} Inflation process: To construct $\G_{n+1}'$, take $d$ copies of $\G_n'$ and label each copy by a letter in $X$. The vertex set of the $x$-th copy $\G_n' \times x$,  can be identified with the set $X^n \times \{x\}.$ Two vertices $(u,x)$ and $(v,y)$ are connected by an edge if and only if $u, v \in \pc$ and $\left\{(u, x), (v, y)\right\} \in E_M.$
\end{mythm}

To construct the Schreier graph $\G_{n+1}$ from the Tile graph $\G_{n+1}'$, we add additional edges between post-critical vertices. Specifically, connect post-critical vertices $p$ and $q$ in $\G_{n+1}'$ by an edge whenever there exists $s \in S$ such that $s(p) = q$ and $s|_p \ne \id$. Note that in such cases, the section $s|_p$ belongs to $S'$, where $S'$ is the set of reachable states.

In the Definition~\ref{def:model} and Theorem~\ref{Theorem:bon-inflation} one can observe the role of the alphabet $X$. After taking $|X|$ copies of $\G_r'$, Theorem~\ref{Theorem:bon-inflation} produces a next level Tile graph $\G_{r+1}'$ and in order to do this, the model graph with vertex set $\pc\times X$ has been used. That is, to construct the Tile graph $\G_{r+n}', (n > 1)$ while starting with the Tile graph $\G_r'$, apply Theorem~\ref{Theorem:bon-inflation} iteratively $n$-times and produce the Tile graph $\G_{r+n}'$. Notice that applying Theorem  \ref{Theorem:bon-inflation} $n$-times iteratively is equivalent to taking $|X^n|$ copies of $\G_r'$ and using edges of the model graph $M_n$. One can connect these copies by identifying model vertices from the $|X^n|$ copies. Let us write this explicitly. For $n-$th iterated inflation graphs we follow Section 1.1 on page no. 94 of \cite{bondarenko2007groups} and Section 3.10 on page no. 110 of \cite{nekrashevych2005self}. More information on Schreier graphs of the Basilica group can be found in \cite{D'angeli2010167}. 
\begin{mydef} \label{def:nth-model} We denote $M_n$ as $n-$th iterated model graph associated to the bounded automaton $\A$ whose vertex set is $V_n = \mathcal{P_A} \times X^n$ and the edge set $E_n$ is given below:
\end{mydef}
\begin{eqnarray}
\label{Equation:model-edges}
E_{M_n} = \left\lbrace \begin{array}{l | l}
& \textnormal{ there exists a path } \cdots e_2e_1 \textnormal{ in the Moore } \\
\{(p,x_n \cdots x_1), & \textnormal{ diagram of } \A, \textnormal{ which ends in the trivial }\\
(q,y_n \cdots y_1)\} & \textnormal{ state and which is labelled by }\\
& px_n \cdots x_1|qy_n \cdots y_1, \textnormal{ where } ~~p,q \in \mathcal{P_A} \textnormal{ and } \\
&~~x_n \cdots x_1, ~y_n \cdots y_1 \in X^n \end{array} \right\rbrace.
\end{eqnarray}
Now we describe the $n$th iterated {\it inflation} in the Corollary given below.
\begin{mycor}\label{Corollary:iterated-inflation} To construct the Tile graph $\G_{r+n}'$ take $|X|^n$ copies of the Tile graph $\G_r'$, identify their sets of vertices with $X^r \times x_n \cdots x_1$ for all $ x_n \cdots x_1 \in X^n$  and connect two vertices $(u, x_n \cdots x_1)$ and $(v, y_n \cdots y_1)$ by an edge if $u, v \in \mathcal{P_A}$ and $\{(u,  x_n \cdots x_1); (v,  y_n \cdots y_1)\} \in E_{M_n}.$
\end{mycor}
\begin{proof}
	It follows immediately by applying iteratively Theorem \ref{Theorem:bon-inflation} $n$-times.
\end{proof}
\subsection{Examples}
\begin{figure}[!htb]
	\centering
    \begin{subfigure}[b]{0.45\textwidth}
        \begin{tikzpicture}[scale=.8]
		\tikzstyle{every state}=[draw=black, thick,fill=none]
		\node[state] at (-2,2)  (a) {$a$};
		\node[state] at (-2,-2) (b) {$b$};
		\node[state] at (2,0) (id) {$\id$};
		\draw[thick,->] (a)  to  (id) ;
		\node[above] at (0,1) {$1|1$};
		\draw[thick,->] (b)  to  (id) ;
		\node[below] at (0,-1){$1|0$};
		\draw[thick,->] (a)  to [out= 240, in=120] (b) ;
		\node[] at (-.75,0) {$0|1$};
		\draw[thick,->] (b)  to [out= 60, in=300] (a) ;
		\node[] at (-2.3,0) {$0|0$};
		\path[->] 	(id) edge [loop right] node[right] {$0|0, 1|1$} ();
		\end{tikzpicture}
        \caption{Generating automaton of $\mathcal{B}$}
		\label{Figure:basilica-aut:left}
    \end{subfigure}
    \begin{subfigure}[b]{0.45\textwidth}
        \begin{tikzpicture}[scale=.7]
        \draw[dashed, -] (-1,0) -- (1,0) -- (2,1.5);
		\draw[fill] (-1,0) circle [radius=0.035] node[below]  {$p_30$};
		\draw[fill] (-3,0) circle [radius=0.035] node[above]  {$\underline{p_20 = p_3}$};
		\draw[fill] (-2,1.5) circle [radius=0.035] node[above]  {$\underline{p_10 = p_1}$};

        \draw[fill] (1,0) circle [radius=0.035] node[below]  {$p_11$};
		\draw[fill] (3,0) circle [radius=0.035] node[right]  {$\underline{p_31 = p_2}$};
		\draw[fill] (2,1.5) circle [radius=0.035] node[right]  {$p_21$};
		\draw (-3,0) circle [radius=0.1];
		\draw (-2,1.5) circle [radius=0.1];
		\draw (3,0) circle [radius=0.1];
		\end{tikzpicture}
        \caption{Model graph of $\mathcal{B}$}
		\label{Figure:basilica-model:right}
    \end{subfigure}
	\caption{The Basilica group $\mathcal{B}$ with post-critical set $\mathcal{P_B} = \bigl\{ p_1 = 0^{-\omega}, p_2 = (01)^{-\omega}, p_3 = (10)^{-\omega}\bigr\}$}
	\label{Figure:basilica}
\end{figure}
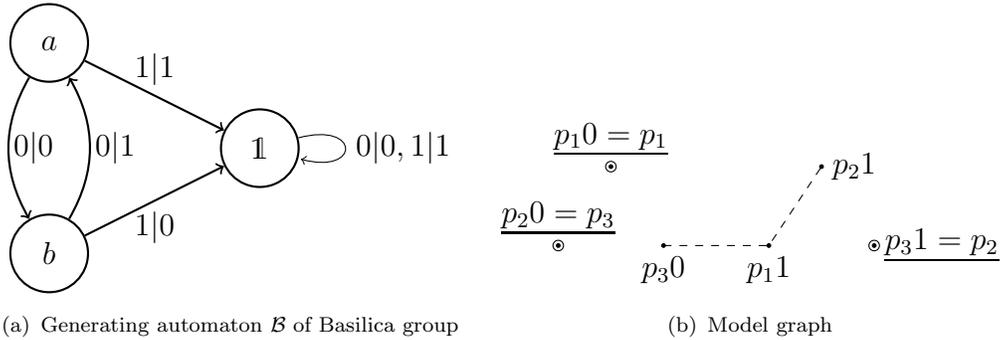
Here we shall apply the inflation to the Schreier graphs of the Basilica group. 
\begin{myexa} Basilica group $\B$\label{basilica-ex} \end{myexa} 
This group is generated by three state automaton as shown in Figure  \ref{Figure:basilica-aut:left}. We use the notation $\G_n'^{\B}$ and $\G_n^{\B}$ for the $n$--level Tile and Schreier graph of $\B$, respectively. The corresponding {\it model graph} is shown in Figure \ref{Figure:basilica-model:right}. The Figure \ref{Figure:Basilicatile123} explains inflation. The dashed edges shown in Figure \ref{Figure:basilica-model:right} and \ref{Figure:Basilicatile123} are model edges. In Figure~\ref{Figure:basilica-model:right}, the underlined vertices correspond to post-critical sequences, while in $\G_3'^{\B}$ of Figure~\ref{Figure:BasilicaSchreier123} they correspond to post-critical vertices. In order to construct the Schreier graph $\G_n^{\B}$ we take the Tile graph $\G_n'^{\B}$ and we add edges (see the dark edges in Figure \ref{Figure:BasilicaSchreier123}) between post-critical vertices $u$ and $v$ if $s(u) = v$ for some $s \in S$ with $s|_u \neq \id$. 
\begin{figure}[!htb]
	\begin{tikzpicture}[scale=1.6]
	\draw[thin,-] (1,1)  to [out=225,in=315] (0,1) ;
	\draw[thin] (1.1,1) circle [radius=0.1];
	\draw[fill] (0,1) circle [radius=0.035] node[above left]  {$p_1=p_3=0$};
	\draw[fill] (1,1) circle [radius=0.035] node[above]  {$1=p_2$};
	\draw[-] (0.5,0.5) node[]{$\G_1'^{\B}$}; 
	%..........................................................................................
	%($\G_2$)
	\draw[thin,-] (2,1) to [out=45,in=135] (3,1);
	\draw[dashed,-] (4,1) to [out=225,in=315] (3,1) ;
	\draw[dashed,-] (4,1)  to [out=45,in=135] (5,1) ;
	\draw[thin,-] (5,1)  to [out=225,in=315] (4,1) ;
	
	\draw[thin] (5.1,1) circle [radius=0.1] ;
	\draw[thin] (1.9,1) circle [radius=0.1] ;
	\draw[fill] (2,1) circle [radius=0.035] node[below=4pt]  {$10$};
	\draw[fill] (3,1) circle [radius=0.035] node[above=5pt]  {$p_30=00$};
	\draw[fill] (4,1) circle [radius=0.035] node[above right=5pt]  {$01=p_11$};
	\draw[fill] (5,1) circle [radius=0.035] node[below=7pt]  {$p_21=11$};
	\draw[-] (3.5,0.5) node[]{$\G_2'^{\B}$}; 
	%........................................................................................
	%($\G_3$)
	\draw[thin,-] (0.4,-1) to [out=45,in=135] (1.2,-1)  ;
	\draw[thin,-] (1.2,-1) to [out=225,in=315] (0.4,-1) ;
	\draw[thin,-] (1.2,-1) to [out=45,in=135] (2,-1) ;

	\draw[thin,-] (2,-1) to [out=45,in=225] (2.5,-0.5) ;         
	\draw[dashed,-] (2.5,-0.5) to [out=-45,in=135] (3,-1) ;
	\draw[thin,-] (3,-1) to [out=225,in=45] (2.5,-1.5) ;

	\draw[dashed,-] (3,-1) to [out=45,in=135] (3.8,-1) ;
	\draw[thin,-] (3.8,-1) to [out=225,in=315] (3,-1) ;
	\draw[thin,-] (3.8,-1) to [out=45,in=135] (4.6,-1) ;
	\draw[thin,-] (4.6,-1) to [out=225,in=315] (3.8,-1) ;

	\draw[thin] (4.7,-1) circle [radius=0.1] ;
	\draw[thin] (0.3,-1) circle [radius=0.1] ;
	\draw[thin] (2.5,-0.4) circle [radius=0.1] ;
	\draw[thin] (2.5,-1.6) circle [radius=0.1] ;
	
	\draw[fill] (.4,-1) circle [radius=0.035] node[below=6pt]  {$110$};
	\draw[fill] (1.2,-1) circle [radius=0.035] node[below=6pt]  {$\underline{010}$};
	\draw[fill] (2,-1) circle [radius=0.035] node[below=6pt]  {$\underline{000}$};
	\draw[fill] (2.5,-0.5) circle [radius=0.035] node[left=4pt]  {$p_30=100$};
	\draw[fill] (2.5,-1.5) circle [radius=0.035] node[below=6pt]  {$\underline{101}$};
	\draw[fill] (3,-1) circle [radius=0.035];
	\draw[-] (3.3,-1) node[below=7pt]  {$001=p_11$};
	\draw[fill] (3.8,-1) circle [radius=0.035] node[above=7pt]  {$p_21=011$};
	\draw[fill] (4.6,-1) circle [radius=0.035] node[below=6pt]  {$111$};
	\draw[-] (2.5,-2.5) node[]{$\G_3'^{\B}$}; 
	%.......................................................................................
	\end{tikzpicture}
	\centering
    \caption{\small The graphs $\G_1'^{\B}$, $\G_2'^{\B}$, and $\G_3'^{\B}$ are the Tile graphs of $\mathcal{B}$ over $X$, $X^2$, and $X^3$, respectively. The corresponding edge set is $\textstyle E_M^B = \biggl\{\{p_{11}, p_{30}\}, \{p_{11}, p_{21}\}\biggr\}$.}

	\label{Figure:Basilicatile123}
\end{figure}
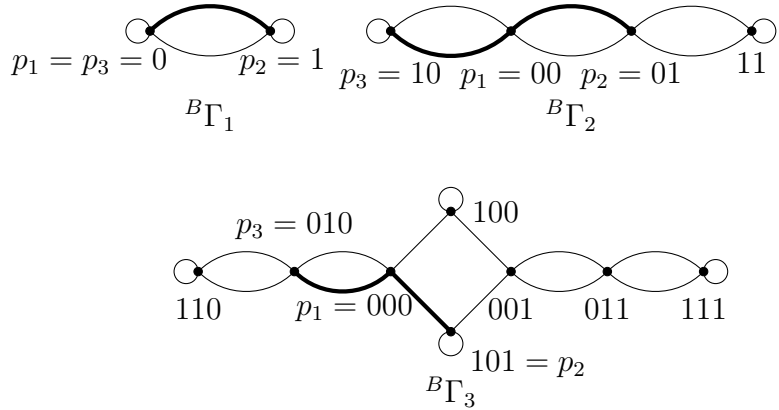
\begin{figure}[!htb]
	\begin{tikzpicture}[scale=1.6]
	%($\G_1$)
	\draw[ultra thick,-] (0,1)  to [out=45,in=135] (1,1) ;
	\draw[thin,-] (1,1)  to [out=225,in=315] (0,1) ;
	\draw[thin] (1.1,1) circle [radius=0.1] ;
	\draw[thin] (-0.1,1) circle [radius=0.1] ;
	\draw[fill] (0,1) circle [radius=0.035]; 
	\draw[-] (-0.5,0.5) node[above=2pt]  {$p_1=p_3=0$};
	\draw[fill] (1,1) circle [radius=0.035]; 
	\draw[-] (1.1,0.5) node[above=2pt]  {$p_2=1$};
	\draw[-] (0.5,0.3) node[]{$\G_1^{\B}$}; 
	%..........................................................................................
	%($\G_2$)
	\draw[thin,-] (2,1)  to [out=45,in=135] (3,1) ;
	\draw[ultra thick,-] (3,1)  to [out=225,in=315] (2,1) ;
	\draw[ultra thick,-] (3,1)  to [out=45,in=135] (4,1) ;
	\draw[thin,-] (4,1)  to [out=225,in=315] (3,1) ;
	\draw[thin,-] (4,1)  to [out=45,in=135] (5,1) ;
	\draw[thin,-] (5,1)  to [out=225,in=315] (4,1) ;
	
	\draw[thin] (5.1,1) circle [radius=0.1] ;
	\draw[thin] (1.9,1) circle [radius=0.1] ;
	\draw[fill] (2,1) circle [radius=0.035] node[below=8pt]  {$p_3=10$};
	\draw[fill] (3,1) circle [radius=0.035] node[below=8pt]  {$p_1=00$};
	\draw[fill] (4,1) circle [radius=0.035] node[below=8pt]  {$p_2=01$};
	\draw[fill] (5,1) circle [radius=0.035] node[below=4pt]  {$11$};
	\draw[-] (3.5,0.3) node[]{$\G_2^{\B}$}; 
	%........................................................................................
	%($\G_3$)
	\draw[thin,-,-] (0.4,-1)  to [out=45,in=135] (1.2,-1)  ;
	\draw[thin,-,-] (1.2,-1)  to [out=225,in=315] (0.4,-1) ;
	\draw[thin,-] (1.2,-1)  to [out=45,in=135] (2,-1)  ;
	\draw[ultra thick,-] (2,-1)  to [out=225,in=315] (1.2,-1) ;
	
	\draw[thin,-] (2,-1)  to [out=45,in=225] (2.5,-0.5) ;         
	\draw[thin,-] (2.5,-0.5)  to [out=-45,in=135] (3,-1) ;
	\draw[thin,-] (3,-1)  to [out=225,in=45] (2.5,-1.5) ; 
	\draw[ultra thick,-] (2.5,-1.5)  to [out=135,in=-45] (2,-1) ;
	
	\draw[thin,-] (3,-1)  to [out=45,in=135] (3.8,-1) ;
	\draw[thin,-] (3.8,-1)  to [out=225,in=315] (3,-1) ;
	\draw[thin,-] (3.8,-1)  to [out=45,in=135] (4.6,-1) ;
	\draw[thin,-] (4.6,-1)  to [out=225,in=315] (3.8,-1) ;

	\draw[thin] (4.7,-1)  circle [radius=0.1];
	\draw[thin] (0.3,-1)  circle [radius=0.1] ;
	\draw[thin] (2.5,-0.4)  circle [radius=0.1] ;
	\draw[thin] (2.5,-1.6) circle [radius=0.1];
	
	\draw[fill] (.4,-1) circle [radius=0.035] node[below=6pt]  {$110$};
	\draw[fill] (1.2,-1) circle [radius=0.035] node[above=8pt]  {$p_3 =010$};
	\draw[fill] (2,-1) circle [radius=0.035];
	\draw[-] (1.7,-1.3) node[]  {$p_1=000$};
	\draw[fill] (2.5,-0.5) circle [radius=0.035] node[right=4pt]  {$100$};
	\draw[fill] (2.5,-1.5) circle [radius=0.035] node[below right=3pt]  {$101= p_2$};
	\draw[fill] (3,-1) circle [radius=0.035] node[below=6pt]  {$001$};
	\draw[fill] (3.8,-1) circle [radius=0.035] node[below=6pt]  {$011$};
	\draw[fill] (4.6,-1) circle [radius=0.035] node[below=6pt]  {$111$};
	\draw[-] (2.5,-2) node[]{$\G_3^{\B}$}; 
	%.......................................................................................
	\end{tikzpicture}
	\centering
	\caption{The graphs $\G_1^{\B},\G_2^{\B}$ and $\G_3^{\B}$ are the Schreier graphs of the Basilica group over $X, X^2$ and $X^3$ respectively.}
	\label{Figure:BasilicaSchreier123}
\end{figure}

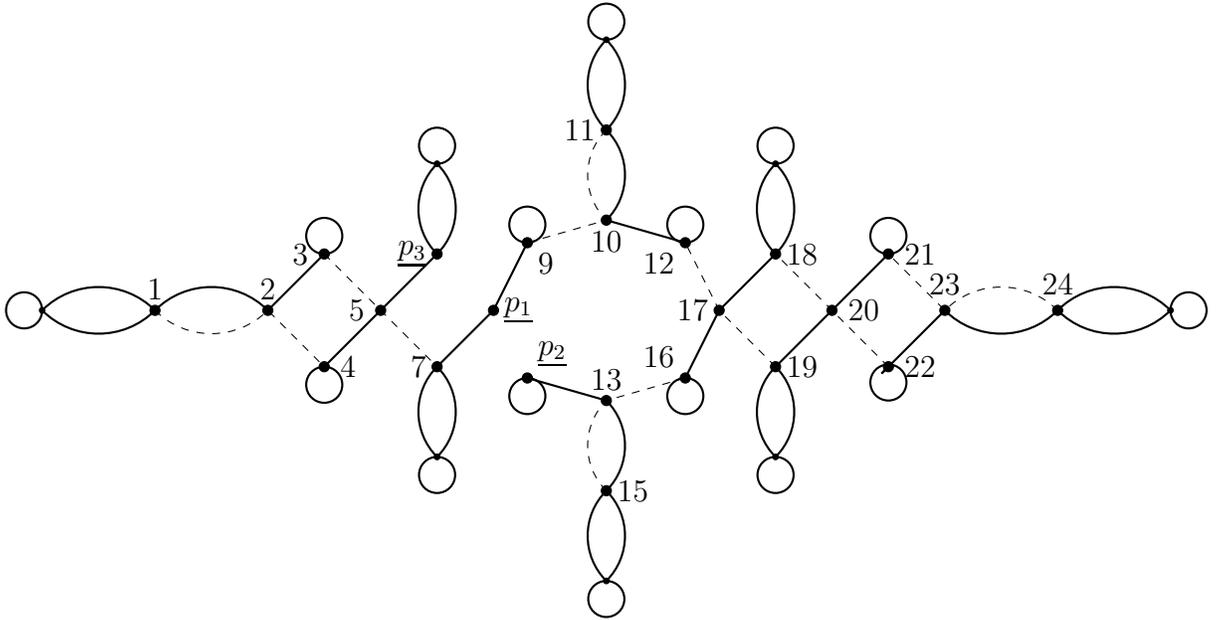
\begin{figure}[!htb]
	\centering
	\begin{tikzpicture}[scale=1.25]
	%.........$\G_5$.....................................................
	%..vertices.............................................................. 
	\draw[fill] (0,1) circle [radius=0.025];
	\draw[fill] (1,1) circle [radius=0.045] node[above] {1};
	\draw[fill] (2,1) circle [radius=0.045] node[above] {2};
	\draw[fill] (2.5,1.5) circle [radius=0.045] node[left=2pt] {3};
	
	\draw[fill] (2.5,0.5) circle [radius=0.045] node[right=2pt] {4};
	\draw[fill] (3,1) circle [radius=0.045] node[left=2pt] {5};
	\draw[fill] (3.5,1.5) circle [radius=0.045] node[left] {$\underline{p_3}$};
	\draw[fill] (3.5,2.3) circle [radius=0.025];
	
	\draw[fill] (3.5,0.5) circle [radius=0.045] node[left] {7};
	\draw[fill] (3.5,-0.3) circle [radius=0.025];
	\draw[fill] (4,1) circle [radius=0.045] node[right] {$\underline{p_1}$};
	\draw[fill] (4.3,1.6) circle [radius=0.045] node[below right] {9};
	
	\draw[fill] (5,1.8) circle [radius=0.045] node[below] {10};
	\draw[fill] (5,2.6) circle [radius=0.045] node[left] {11};
	\draw[fill] (5,3.4) circle [radius=0.025];
	\draw[fill] (5.7,1.6) circle [radius=0.045] node[below left] {12};
	
	\draw[fill] (5,0.2) circle [radius=0.045] node[above] {13};
	\draw[fill] (4.3,0.4) circle [radius=0.045]node[above right] {$\underline{p_2}$};
	\draw[fill] (5,-0.6) circle [radius=0.045] node[right] {15};
	\draw[fill] (5,-1.4) circle [radius=0.025];
	
	\draw[fill] (5.7,0.4) circle [radius=0.045] node[above left] {16};
	\draw[fill] (6,1) circle [radius=0.045] node[left] {17};
	\draw[fill] (6.5,1.5) circle [radius=0.045] node[right] {18};
	\draw[fill] (6.5,2.3) circle [radius=0.025];
	
	\draw[fill] (6.5,0.5) circle [radius=0.045] node[right] {19};
	\draw[fill] (6.5,-0.3) circle [radius=0.025];
	\draw[fill] (7,1) circle [radius=0.045] node[right=2pt] {20};
	\draw[fill] (7.5,1.5) circle [radius=0.045] node[right=2pt] {21};
	
	\draw[fill] (7.5,0.5) circle [radius=0.045] node[right=2pt] {22};
	\draw[fill] (8,1) circle [radius=0.045] node[above=2pt] {23};
	\draw[fill] (9,1) circle [radius=0.045] node[above=2pt] {24};
	\draw[fill] (10,1) circle [radius=0.025];
	%........................Edges...............................................................
	
	%Gamma_2, 1st sheet................................................
	\draw[thick,-] (0,1) to [out=45,in=135] (1,1) to [out=45,in=135](2,1) to [out = 45 , in =225 ] (2.5,1.5);
	\draw[thick,-] (0,1)  to [out=-45,in=225] (1,1)  ;
	%Gamma_2, 2nd sheet......................................................
	\draw[thick,-](2.5,0.5) to [out = 45, in =225 ](3,1) to [out = 45 , in = 225 ](3.5,1.5) to [out = 45 , in = -45](3.5,2.3) to [out =225 , in = 135](3.5,1.5); 
	%Gamma_2, third sheet......................................................
	\draw[thick,-](4.3,1.6) to (4,1) to (3.5,0.5) to [out =225 , in = 135](3.5,-0.3) to  [out =45 , in =-45](3.5,0.5) ; 
	%Gamma_2, 4th sheet......................................................
	\draw[thick,-](5.7,1.6) to (5,1.8);
	\draw[thick,-](5,2.6) to [out =135 , in =225](5,3.4) to [out =-45 , in = 45](5,2.6) to [out =-45 , in = 45](5,1.8) ; 
	%Gamma_2, 5th sheet......................................................
	\draw[thick,-](4.3,0.4) to (5,0.2);
	\draw[thick,-](5,-0.6) to  [out =225 , in =135](5,-1.4) to [out =45 , in = -45](5,-0.6) to [out =45 , in = -45](5,0.2) ; 
	%Gamma_2, 6th sheet......................................................
	\draw[thick,-](5.7,0.4) to (6,1) to (6.5,1.5) to [out =135 , in = 225](6.5,2.3) to  [out =-45 , in =45](6.5,1.5) ; 
	%Gamma_2, 7th sheet......................................................
	\draw[thick,-](7.5,1.5) to [out = 225, in =45 ](7,1) to [out =225 , in =45 ](6.5,0.5) to [out =225 , in =135](6.5,-0.3) to [out =45 , in =-45](6.5,0.5); 
	%Gamma_2, 8th sheet......................................................
	\draw[thick,-]  (9,1) to [out=45,in=135](10,1);
	\draw[thick,-] (10,1) to [out=225,in=-45] (9,1) to [out=225,in=-45](8,1) to [out = 225 , in =-135 ]  (7.5,0.5);
	%loops.................................................................
	\draw[thick] (-.16,1) circle [radius=0.16];
	\draw[thick] (2.5,1.66) circle [radius=0.16];
	\draw[thick] (2.5,0.34) circle [radius=0.16];
	\draw[thick] (3.5,2.46) circle [radius=0.16];
	\draw[thick] (3.5,-0.46) circle [radius=0.16];
	\draw[thick] (4.3,1.76) circle [radius=0.16];
	\draw[thick] (4.3,0.24) circle [radius=0.16];
	\draw[thick] (5.7,1.76) circle [radius=0.16];
	\draw[thick] (5.7,0.24) circle [radius=0.16];
	\draw[thick] (5,3.56) circle [radius=0.16];
	\draw[thick] (5,-1.56) circle [radius=0.16];
	\draw[thick] (6.5,2.46) circle [radius=0.16];
	\draw[thick] (6.5,-0.46) circle [radius=0.16];
	\draw[thick] (7.5,1.66) circle [radius=0.16];
	\draw[thick] (7.5,0.36) circle [radius=0.16];
	\draw[thick] (10.16,1) circle [radius=0.16];
	%..............................Lifts...........................................
	\draw[dashed,-] (1,1) to [out=-45 , in= 225](2,1) to (2.5,0.5);
	\draw[dashed,-] (2.5,1.5) to (3,1) to (3.5,0.5);
	\draw[dashed,-] (4.3,1.6) to (5,1.8) to [out=135 ,in=225 ](5,2.6);
	\draw[dashed,-] (5.7,0.4) to (5,0.2) to [out=225 ,in=135 ](5,-0.6);
	\draw[dashed,-] (6.5,0.5) to (6,1) to (5.7,1.6);
	\draw[dashed,-] (7.5,0.5) to (7,1) to (6.5,1.5);
	\draw[dashed,-] (9,1) to [out=135 , in= 45](8,1) to (7.5,1.5);
	%..............................................................................
	\end{tikzpicture}
	\vspace{2cm}\\
	\caption{The Tile graph $\G_5'^{\B}$ of Basilica group over $X^5$, $\mathcal{P_B} = \{ p_1 = (0)^{-\omega}, p_2 = (01)^{-\omega}, p_3 = (10)^{-\omega}\}.$}
	\label{Figure:B-tile5}
\end{figure}

\begin{figure}[!htb]
	\centering
	\begin{tikzpicture}[scale=1.25]
	%.........$\G_5$.....................................................
	%..vertices.............................................................. 
	\draw[fill] (0,1) circle [radius=0.025];
	\draw[fill] (1,1) circle [radius=0.025];
	\draw[fill] (2,1) circle [radius=0.025];
	\draw[fill] (2.5,1.5) circle [radius=0.025];
	\draw[fill] (2.5,0.5) circle [radius=0.025];
	\draw[fill] (3,1) circle [radius=0.025];
	\draw[fill] (3.5,1.5) circle [radius=0.045] node[left] {$\underline{p_3}$};
	\draw[fill] (3.5,2.3) circle [radius=0.025];
	\draw[fill] (3.5,0.5) circle [radius=0.025];
	\draw[fill] (3.5,-0.3) circle [radius=0.025];
	\draw[fill] (4,1) circle [radius=0.045] node[right] {$\underline{p_1}$};
	\draw[fill] (4.3,1.6) circle [radius=0.025];
	\draw[fill] (5,1.8) circle [radius=0.025];
	\draw[fill] (5,2.6) circle [radius=0.025];
	\draw[fill] (5,3.4) circle [radius=0.025];
	\draw[fill] (5.7,1.6) circle [radius=0.025];
	\draw[fill] (5,0.2) circle [radius=0.025];
	\draw[fill] (4.3,0.4) circle [radius=0.045]node[above right] {$\underline{p_2}$};
	\draw[fill] (5,-0.6) circle [radius=0.025];
	\draw[fill] (5,-1.4) circle [radius=0.025];
	\draw[fill] (5.7,0.4) circle [radius=0.025];
	\draw[fill] (6,1) circle [radius=0.025];
	\draw[fill] (6.5,1.5) circle [radius=0.025];
	\draw[fill] (6.5,2.3) circle [radius=0.025];
	\draw[fill] (6.5,0.5) circle [radius=0.025];
	\draw[fill] (6.5,-0.3) circle [radius=0.025];
	\draw[fill] (7,1) circle [radius=0.025];
	\draw[fill] (7.5,1.5) circle [radius=0.025];
	\draw[fill] (7.5,0.5) circle [radius=0.025];
	\draw[fill] (8,1) circle [radius=0.025];
	\draw[fill] (9,1) circle [radius=0.025];
	\draw[fill] (10,1) circle [radius=0.025];
	%........................Edges...............................................................
	
	%Gamma_2, 1st sheet................................................
	\draw[thick,-] (0,1) to [out=45,in=135] (1,1) to [out=45,in=135](2,1) to [out = 45 , in =225 ] (2.5,1.5);
	\draw[thick,-] (0,1)  to [out=-45,in=225] (1,1)  ;
	%Gamma_2, 2nd sheet......................................................
	\draw[thick,-](2.5,0.5) to [out = 45, in =225 ](3,1) to [out = 45 , in = 225 ](3.5,1.5) to [out = 45 , in = -45](3.5,2.3) to [out =225 , in = 135](3.5,1.5); 
	%Gamma_2, third sheet......................................................
	\draw[thick,-](4.3,1.6) to (4,1) to (3.5,0.5) to [out =225 , in = 135](3.5,-0.3) to  [out =45 , in =-45](3.5,0.5) ; 
	%Gamma_2, 4th sheet......................................................
	\draw[thick,-](5.7,1.6) to (5,1.8);
	\draw[thick,-](5,2.6) to [out =135 , in =225](5,3.4) to [out =-45 , in = 45](5,2.6) to [out =-45 , in = 45](5,1.8) ; 
	%Gamma_2, 5th sheet......................................................
	\draw[thick,-](4.3,0.4) to (5,0.2);
	\draw[thick,-](5,-0.6) to  [out =225 , in =135](5,-1.4) to [out =45 , in = -45](5,-0.6) to [out =45 , in = -45](5,0.2) ; 
	%Gamma_2, 6th sheet......................................................
	\draw[thick,-](5.7,0.4) to (6,1) to (6.5,1.5) to [out =135 , in = 225](6.5,2.3) to  [out =-45 , in =45](6.5,1.5) ; 
	%Gamma_2, 7th sheet......................................................
	\draw[thick,-](7.5,1.5) to [out = 225, in =45 ](7,1) to [out =225 , in =45 ](6.5,0.5) to [out =225 , in =135](6.5,-0.3) to [out =45 , in =-45](6.5,0.5); 
	%Gamma_2, 8th sheet......................................................
	\draw[thick,-]  (9,1) to [out=45,in=135](10,1);
	\draw[thick,-] (10,1) to [out=225,in=-45] (9,1) to [out=225,in=-45](8,1) to [out = 225 , in =-135 ]  (7.5,0.5);
	%loops.................................................................
	\draw[thick] (-.16,1) circle [radius=0.16];
	\draw[thick] (2.5,1.66) circle [radius=0.16];
	\draw[thick] (2.5,0.34) circle [radius=0.16];
	\draw[thick] (3.5,2.46) circle [radius=0.16];
	\draw[thick] (3.5,-0.46) circle [radius=0.16];
	\draw[thick] (4.3,1.76) circle [radius=0.16];
	\draw[thick] (4.3,0.24) circle [radius=0.16];
	\draw[thick] (5.7,1.76) circle [radius=0.16];
	\draw[thick] (5.7,0.24) circle [radius=0.16];
	\draw[thick] (5,3.56) circle [radius=0.16];
	\draw[thick] (5,-1.56) circle [radius=0.16];
	\draw[thick] (6.5,2.46) circle [radius=0.16];
	\draw[thick] (6.5,-0.46) circle [radius=0.16];
	\draw[thick] (7.5,1.66) circle [radius=0.16];
	\draw[thick] (7.5,0.36) circle [radius=0.16];
	\draw[thick] (10.16,1) circle [radius=0.16];
	%..............................Lifts...........................................
	\draw[thick,-] (1,1) to [out=-45 , in= 225](2,1) to (2.5,0.5);
	\draw[thick,-] (2.5,1.5) to (3,1) to (3.5,0.5);
	\draw[ultra thick,-] (3.5,1.5) to (4,1) to (4.3,0.4);
	\draw[thick,-] (4.3,1.6) to (5,1.8) to [out=135 ,in=225 ](5,2.6);
	\draw[thick,-] (5.7,0.4) to (5,0.2) to [out=225 ,in=135 ](5,-0.6);
	\draw[thick,-] (6.5,0.5) to (6,1) to (5.7,1.6);
	\draw[thick,-] (7.5,0.5) to (7,1) to (6.5,1.5);
	\draw[thick,-] (9,1) to [out=135 , in= 45](8,1) to (7.5,1.5);
	%..............................................................................
	\end{tikzpicture} 
	\caption{The Schreier graph $\G_5^{\B}$ of Basilica group over $X^5$.}
	\label{Figure:B-Sch5}
\end{figure}
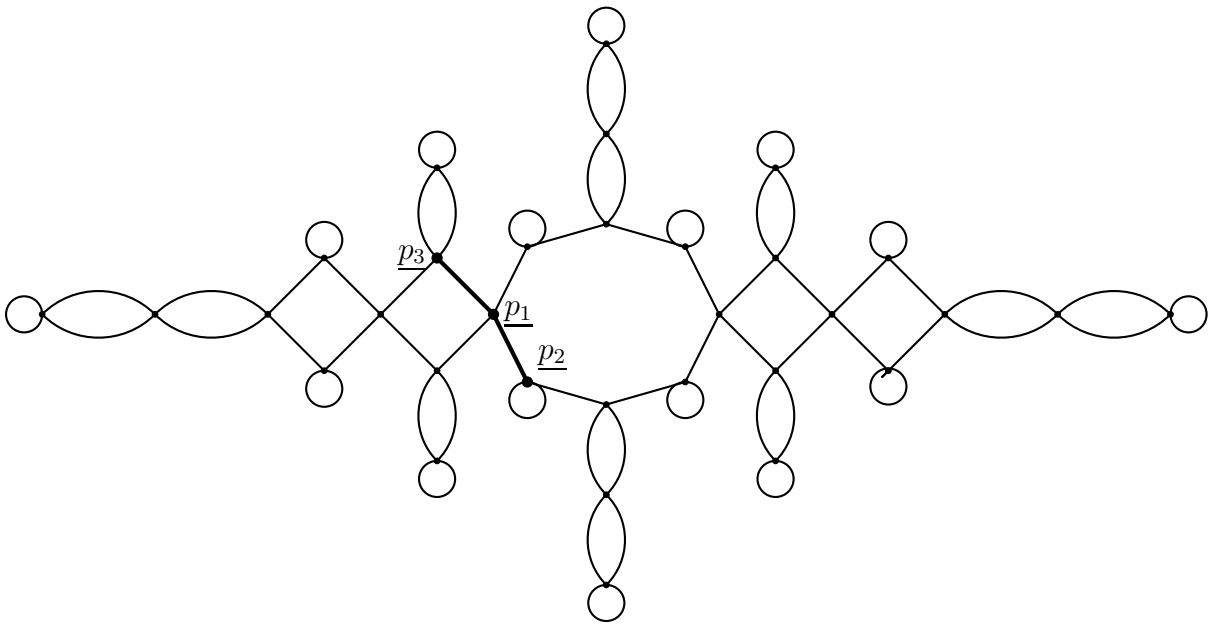

\begin{myexa}  We shall construct the Tile graph $\G_5'^{\B}$ and hence Schreier graph $\G_5^{\B}$ by applying Corollary \ref{Corollary:iterated-inflation}. \end{myexa}
In order to do this we start with the Tile graph $\G_2'^{\B}$. See Figures \ref{Figure:B-tile5} and \ref{Figure:B-Sch5}. The numbered vertices and dashed edges are shown in the Figure  \ref{Figure:B-tile5} are the vertices and the edges of third iterated model graph $M_3$. For the vertex numbering see below Table \ref{Table:B-tile5vertices}.
\begin{table}[!htb]
	\caption{Vertices of $3$rd iterated model graph $M_3$.}
	
	\centering
	{\begin{tabular}{c|c|c|c} 
			Numbering & Vertices  &  Numbering &  Vertices \\
			\hline
			\hline
			$1,2,3$ & $01110, 00110, 10110$ & $13,p_2,15$ & $00101, 10101, 01101$ \\
			$4,5,p_3$ & $ 10010, 00010, 01010$ & $16,17,18$ & $10001, 00001, 01001$ \\
			$7,p_1,9$ & $ 01000, 00000, 10000$ & $19,20,21$ & $01011, 00011, 10011$ \\
			$10,11,12$ & $00100, 01100, 10100 $ & $22,23,24$ & $10111, 00111, 01111$ \\
	\end{tabular}}
	\label{Table:B-tile5vertices}
\end{table}

\section{Generalized Replacement Product}\label{sec:replacement}

The replacement product of graphs is a classical construction in graph theory. Given two regular graphs $G_1$ and $G_2$ of degrees $d_1$ and $d_2$, respectively, their replacement product, denoted $G_1 \textcircled{r} G_2$, is a regular graph of degree $d_2+1$ (see \cite{abdollahi2013one, D'Angeli2015JGT} for detailed treatments).

In the setting of Schreier graphs, an analogue of this construction was first introduced in \cite{shaikh2016zeta} for the Schreier graphs $\G_n^{\B}$ of the Basilica group (see Definition~2.1), where it was used to describe Galois coverings (see Proposition~2.1 therein). Motivated by this analogy, we continue to use the term \emph{generalized replacement product} for our construction, even though the degree behavior differs from the classical case.

In this section we extend the generalized replacement product to Schreier graphs $\G_n$ associated with any group $G$ generated by a bounded automaton $\A$. Our construction preserves $d$-regularity within the Schreier graph framework and naturally reflects the covering structures observed in the Basilica case.

Let $\G = (V, E)$ be a finite, connected, $d$-regular graph that may contain loops and multiple edges. Suppose we assign a set of $d$ colors to the edges, identified with the set of natural numbers $[d] := \{1, 2, \dots, d\}$. For an edge $e = \{v, v'\} \in E$, let the color adjacent to $v$ be denoted by $s$ and the color adjacent to $v'$ by $s'$, with the requirement that distinct edges incident to the same vertex $v$ receive distinct colors near $v$.

Then, the \emph{rotation map} $Rot_{\G} : V \times [d] \to V \times [d]$ is defined by
$$Rot_{\G}(v, s) = (v', s'),$$
for all $v, v' \in V$ and $s, s' \in [d]$. Note that in general we may have $s \ne s'$.

By definition, the composition $Rot_{\G} \circ Rot_{\G}$ is the identity map. Since each edge $e = \{v, v'\}$ is colored by $s$ near $v$ and by $s'$ near $v'$, we refer to $\G$ as a \emph{bi-colored} graph.

Let $\G_n$ and $\G_r$ be two Schreier graphs of a group $G$ generated by a bounded automaton $\A$ with generating set $S$. For an edge $\{v, v'\}$ in these graphs, we say the color near $v$ is $s$ and near $v'$ is $s'$, provided that $o(s, v) = v'$ and $t(s, v) = s'$, where $s' = s^{-1}$ in the group-theoretic sense.

Recall that the notation $s|_u = \id$ indicates that the section of $s$ at $u$ is the trivial state. To define the generalized replacement product, we recall the description of the Tile graph $\G_r'$ and the non-Tile edges of $\G_r$ given in Equations~\eqref{eq:nonSchreier} and \eqref{eqn:non-tile}.

\begin{mydef}  \label{Definition:gen-rep}
The \emph{generalized replacement product} $\G_n \g \G_r$ is a $|S|$-regular graph with vertex set $X^{n+r} = X^n \times X^r$. Its edges are defined via the following rotation map. Let $(v, u) \in X^n \times X^r$ and $s \in S$:
\begin{align}
\label{Equation:gr-rot1}
Rot_{\G_n \g \G_r}((v, u), s) &= ((v, s(u)), s^{-1}) & \textnormal{if } s|_u = \id, \\
\label{Equation:gr-rot2}
Rot_{\G_n \g \G_r}((v, u), s) &= ((s|_u(v), s(u)), s^{-1}) & \textnormal{if } s|_u \ne \id \textnormal{ and } s|_{uv} = \id, \\
\label{Equation:gr-rot3}
Rot_{\G_n \g \G_r}((v, u), s) &= ((s|_u(v), s(u)), s^{-1}) & \textnormal{if } s|_u \ne \id \textnormal{ and } s|_{uv} \ne \id.
\end{align} 
\end{mydef}

The vertex set $X^{n+r}$ of $\G_n \g \G_r$ can be viewed as partitioned into \emph{sheets}, each indexed by a vertex $v \in X^n$ of $\G_n$. The $v$-th sheet consists of the set $\{(v, u) : u \in X^r\}$. Within each sheet, we place a copy of the Tile graph $\G_r'$ and retain its internal edges unchanged. These intra-sheet edges are governed by Equation~\eqref{Equation:gr-rot1} and are referred to as \emph{$v$-sheet edges}.

Edges that connect different sheets are defined by Equations~\eqref{Equation:gr-rot2} and \eqref{Equation:gr-rot3}. We refer to those arising from~\eqref{Equation:gr-rot2} as \emph{model edges}, and those from~\eqref{Equation:gr-rot3} as \emph{Schreier edges}.

The connectedness of both $\G_n$ and $\G_r'$ ensures that the graph $\G_n \g \G_r$ is connected. Proposition~\ref{Proposition:coveringmap} establishes that $\G_n \g \G_r$ is isomorphic to the Schreier graph $\G_{n+r}$ of the group $G$.

It was shown in Proposition~8.1 of \cite{grigorchuk2002torsion} that $\G_{r+1}$ serves as a covering of $\G_r$. We extend this result by showing that the covering is unramified. Furthermore, under the condition~\eqref{item:root} of Theorem~\ref{Theorem:Galois}, the covering is normal, and the corresponding Galois group is cyclic.

\begin{mypro} \label{Proposition:coveringmap}
Let $n, r \geq 1$. Then the following statements hold:
\begin{enumerate}
    \item The graphs $\G_n \g \G_r$ and $\G_{n+r}$ are isomorphic.
    \item The graph $\G_{n+r}$ is an unramified $d^n$-sheeted covering of $\G_r$.
\end{enumerate}\end{mypro}

\begin{proof}
\begin{enumerate}
    \item We define a map $f : \G_n \g \G_r \to \G_{n+r}$ by $f(v, u) = uv$, where $v \in X^n$ and $u \in X^r$. The map is bijective by construction. To show that $f$ preserves adjacency, recall that $\G_n \g \G_r$ contains three types of edges:
    \begin{itemize}
        \item \emph{Sheet edges} (Equation~\eqref{Equation:gr-rot1}): for $s \in S$ with $s|_u = \id$, we have
       $$ Rot_{\G_n \g \G_r}((v, u), s) = ((v, s(u)), s^{-1}).$$
        Since $s|_u = \id$, it follows that $$Rot_{\G_{n+r}}(uv, s) = (s(u)v, s^{-1}),$$ 
        which confirms that $uv$ is adjacent to $s(u)v$ in $\G_{n+r}$, matching the image of $f$.
        \item \emph{Model and Schreier edges} (Equations~\eqref{Equation:gr-rot2} and~\eqref{Equation:gr-rot3}) are treated analogously. In both cases, $Rot_{\G_n \g \G_r}((v, u), s)$ maps to a vertex whose image under $f$ is adjacent to $uv$ in $\G_{n+r}$.
    \end{itemize}
    Hence, $f$ is a graph isomorphism.
    \item Since $\G_n \g \G_r \cong \G_{n+r}$ and $\G_n \g \G_r$ consists of $d^n$ disjoint sheets (indexed by $v \in X^n$), each isomorphic to $\G_r'$, we define a map
    $$ \pi : \G_{n+r} \to \G_r, \quad \pi(uv) = u.$$
    This map is well-defined and surjective on vertices. It is straightforward to check that $\pi$ maps adjacent vertices in $\G_{n+r}$ to adjacent vertices in $\G_r$, and that for each edge in $\G_r$ and each lift point $u \in X^r$, there are exactly $d^n$ distinct preimages in $\G_{n+r}$. Thus, $\pi$ is an unramified covering map. This completes the proof.
\end{enumerate}
\end{proof}	

In Definition~\ref{Definition:gen-rep}, assuming $n = 1$ leads to the following results:

\begin{mypro} \begin{enumerate}
    \item \label{Proposition:r1} The rotation map \eqref{Equation:gr-rot1} yields the disjoint union \\$\displaystyle \biguplus_{x \in X} (x \times \G_r')$, representing $d$ disjoint copies of the Tile graph $\G_r'$, each indexed by $x \in X$.
    
    \item \label{Proposition:r12} Including the rotation map \eqref{Equation:gr-rot2} adds edges between these copies, resulting in the construction of the Tile graph $\G_{r+1}'$.
    
    \item \label{corollary:r123} Furthermore, incorporating the rotation map \eqref{Equation:gr-rot3} connects post-critical vertices of $\G_{r+1}'$, yielding the Schreier graph $\G_{r+1}$.
\end{enumerate} \end{mypro}

\begin{proof}
\begin{enumerate}
    \item Let $(x,u) \in X \times X^r$. From the rotation map \eqref{Equation:gr-rot1}, observe that the first coordinate $x$ remains unchanged under the action of $s \in S$ whenever $s|_u = \id$. Thus, the transformation acts only on the second component $u$ and corresponds exactly to the Tile graph $\G_r'$. Formally, we have a rotation map 
    \[
    Rot_{\G_r}(u, s) = (s(u), s^{-1}), \quad \text{for } s \in S \text{ with } s|_u = \id,
    \]
    which defines the edge set of $\G_r'$. Therefore, each $x \in X$ determines a distinct copy of $\G_r'$ in the product $\G_1 \g \G_r$.
    
    \item When $s|_u \neq \id$ and $s|_{ux} = \id$, the rotation map \eqref{Equation:gr-rot2} creates edges between different sheets indexed by $x \in X$. These edges correspond to transitions in the model graph $M$ (as in Definition~\ref{def:model}) via labels of the form $ux | u'x'$. Thus, by the definition of the inflation process (Theorem~\ref{Theorem:bon-inflation}), the resulting graph is exactly the Tile graph $\G_{r+1}'$.
    
    \item This is a special case of Proposition~\ref{Proposition:coveringmap}, part (1), where applying all three rotation maps reconstructs the full Schreier graph $\G_{r+1}$ from $\G_1 \g \G_r$.
\end{enumerate}
\end{proof}

Recall that Theorem~\ref{Theorem:bon-inflation} describes how to construct the Tile graph $\G_{r+1}'$ from $\G_r'$ via the inflation process. The application of the rotation maps \eqref{Equation:gr-rot1} and \eqref{Equation:gr-rot2} to $\G_r'$ corresponds exactly to this construction. More generally, for any $n, r \geq 1$, applying these rotation maps recursively to $\G_r'$ produces the $n$-th iterated inflation as described in Corollary~\ref{Corollary:iterated-inflation}.

\begin{myexa}
See Figures~\ref{Figure:Gri12} and~\ref{Figure:torsion2.4} for the generalized replacement products of the Schreier graphs $\G_n^{Gr}$ and $\G_n^{BSV}$ associated with the Grigorchuk group and the BSV torsion-free group, respectively. In these figures, the solid (continuous) edges represent sheet edges, the dotted edges indicate model edges, and the dashed edges correspond to Schreier edges.
\end{myexa}
\begin{figure}[!htb]
  \centering
  \begin{minipage}{0.45\textwidth}
    \centering
    \begin{tikzpicture}[scale=1.25]
      \draw[thick,-] (-1,-1) to (1,-1);
      \draw[-] (-1,-1) to [out= 45, in=0] (-1,-0.2) to [out= 180, in=135] (-1,-1);
      \draw[-] (-1,-1) to [out= -45, in=0] (-1,-1.8) to [out= -180, in=-135] (-1,-1);
      \draw[-] (-1,-1) to [out= 135, in=90] (-1.8,-1) to [out= 270, in=225] (-1,-1);
      \draw[-] (1,-1) to [out= 45, in=0] (1,-0.2) to [out= 180, in=135] (1,-1);
      \draw[-] (1,-1) to [out= -45, in=0] (1,-1.8) to [out= -180, in=-135] (1,-1);
      \draw[-] (1,-1) to [out= 45, in=90] (1.8,-1) to [out= 270, in=-45] (1,-1);
      \draw[fill] (-1,-1) circle [radius=0.035];
      \draw[-] (-0.7,-1) node[above]  {$0$};
      \draw[fill] (1,-1) circle [radius=0.035];
      \draw[-] (0.7,-1) node[above]  {$1$};
    \end{tikzpicture}
    \par\vspace{1ex}
    $\G_1^{Gr}$
  \end{minipage}\hfill
  \begin{minipage}{0.45\textwidth}
    \centering
    \begin{tikzpicture}[scale=1.25]
      % Tile graph level 1
      \draw[-] (-1,-1) to (1,-1);
      \draw[-] (-1,-1) to [out= 45, in=0] (-1,-0.2) to [out= 180, in=135] (-1,-1);
      \draw[fill] (-1,-1) circle [radius=0.035] node[left]  {$0$};
      \draw[fill] (1,-1) circle [radius=0.035] node[right] {$1$};
    \end{tikzpicture}
    \par\vspace{1ex}
    $\G_1'^{Gr}$
  \end{minipage}
  % Generalized replacement graph
  \begin{tikzpicture}[scale=1.25]
    % Level 2
    \draw[-] (-3,1) to (-1,1);
    \draw[thick,dotted,-] (-3,1) to [out= 45, in=0] (-3,1.8) to [out= 180, in=135] (-3,1);
    \draw[dashed,-] (-3,1) to [out= -45, in=0] (-3,0.2) to [out= -180, in=-135] (-3,1);
    \draw[dashed,-] (-3,1) to [out= 135, in=90] (-3.8,1) to [out= 270, in=225] (-3,1);
    \draw[-] (-1,1) to [out= -45, in=0] (-1,0.2) to [out= -180, in=-135] (-1,1);
    \draw[thick,dotted,-] (-1,1) to [out= 40, in=140] (1,1);
    \draw[thick,dotted,-] (-1,1) to [out= -40, in=-140] (1,1);
    \draw[-] (1,1) to (3,1);
    \draw[dashed,-] (3,1) to [out= 45, in=0] (3,1.8) to [out= 180, in=135] (3,1);
    \draw[dashed,-] (3,1) to [out= -45, in=0] (3,0.2) to [out= -180, in=-135] (3,1);
    \draw[dashed,-] (3,1) to [out= 45, in=90] (3.8,1) to [out= 270, in=-45] (3,1);
    \draw[-] (1,1) to [out= -45, in=0] (1,0.2) to [out= -180, in=-135] (1,1);
    \draw[fill] (-3,1) circle [radius=0.035];
    \draw[-] (-2.5,1) node[above]  {$(0,1)$};
    \draw[fill] (-1,1) circle [radius=0.035] node[ right]  {$(0,0)$};
    \draw[fill] (1,1) circle [radius=0.035] node[ left]  {$(1,0)$};
    \draw[fill] (3,1) circle [radius=0.035];
    \draw[-] (2.5,1) node[above]  {$(1,1)$};
  \end{tikzpicture}
  \par\vspace{1ex}
  $\G_2^{Gr}$
  \caption{The Schreier graph $\G_1^{Gr}$, the Tile graph $\G_1'^{Gr}$ of the Grigorchuk group, and the generalized replacement graph $\G_1^{Gr} \g \G_1^{Gr}$.}
  \label{Figure:Gri12}
\end{figure}

\begin{figure}[!htb]
\centering
  \begin{minipage}{0.3\textwidth}
    \centering
    \begin{tikzpicture}[scale=.75]
      \draw[thick,-] (-1,1) -- (1,1);
      \draw[thick,-] (1,1) -- (1,-1);
      \draw[thick,-] (1,-1) -- (-1,-1);
      \draw[thick,-] (-1,-1) -- (-1,1);
      \draw[thick,-] (-1,1) to[out=45,in=135] (1,1);
      \draw[thick,-] (1,1) to[out=315,in=45] (1,-1);
      \draw[thick,-] (1,-1) to[out=225,in=315] (-1,-1);
      \draw[thick,-] (-1,-1) to[out=135,in=225] (-1,1);
      \draw[fill] (-1,1) circle[radius=0.035] node[left=4pt]  {0};
      \draw[fill] (1,1) circle[radius=0.035] node[right=4pt] {1};
      \draw[fill] (1,-1) circle[radius=0.035] node[right=4pt] {2};
      \draw[fill] (-1,-1) circle[radius=0.035] node[left=4pt] {3};
    \end{tikzpicture}
    $\G_1^{BSV}$
  \end{minipage}\hfill
  \begin{minipage}{0.3\textwidth}
    \centering
    \begin{tikzpicture}[scale=.75]
      \draw[thick,-] (-1,1) -- (1,1);
      \draw[thick,-] (1,1) -- (1,-1);
      \draw[thick,-] (1,-1) -- (-1,-1);
      \draw[thick,-] (-1,1) to[out=45,in=135] (1,1);
      \draw[thick,-] (1,1) to[out=315,in=45] (1,-1);
      \draw[thick,-] (1,-1) to[out=225,in=315] (-1,-1);
      \draw[fill] (-1,1) circle[radius=0.035] node[left=4pt]  {0};
      \draw[fill] (1,1) circle[radius=0.035] node[right=4pt] {1};
      \draw[fill] (1,-1) circle[radius=0.035] node[right=4pt] {2};
      \draw[fill] (-1,-1) circle[radius=0.035] node[left=4pt] {3};
    \end{tikzpicture}
    $\G_1'^{BSV}$
  \end{minipage} 
	\begin{tikzpicture}[scale=1.5]
	%1st sheet
	\draw[thick,-,-] (-1,2)  to (-1.6,2) ;
	\draw[thick,-,-] (-1.6,2) to (-2,1.6) ;
	\draw[thick,-] (-2,1.6)  to (-2,1) ;
	
	\draw[thick,-] (-1,2)  to [out=135,in=45] (-1.6,2) ;
	\draw[thick,-] (-1.6,2) to [out=180,in=90] (-2,1.6) ;
	\draw[thick,-] (-2,1.6) to [out=225, in=135] (-2,1);
	
	\draw[fill] (-1,2) circle [radius=0.035] node[above right]  {$(0,0)$};
	\draw[fill] (-1.6,2) circle [radius=0.035] node[above left]  {$(0,1)$};
	\draw[fill] (-2,1.6) circle [radius=0.035] node[left=4pt]  {$(0,2)$};
	\draw[fill] (-2,1) circle [radius=0.035] node[left=4pt]  {$(0,3)$};
	
	%2 nd sheet
	\draw[thick,-] (-2,0) to (-2,-.6) ;
	\draw[thick,-] (-2,-.6) to (-1.6,-1) ;
	\draw[thick,-] (-1.6,-1) to (-1,-1) ;
	
	\draw[thick,-] (-2,0) to [out=225,in=135] (-2,-.6) ;
	\draw[thick,-] (-2,-.6) to [out=270,in=180] (-1.6,-1) ;
	\draw[thick,-] (-1.6,-1) to [out=315,in=225] (-1,-1) ;
	
	\draw[fill] (-2,0) circle [radius=0.035] node[left]  {$(1,0)$};
	\draw[fill] (-2,-.6) circle [radius=0.035] node[left=4pt]  {$(1,1)$};
	\draw[fill] (-1.6,-1) circle [radius=0.035] node[left=6pt]  {$(1,2)$};
	\draw[fill] (-1,-1) circle [radius=0.035] node[below=4pt]  {$(1,3)$};
	
	%third sheet
	\draw[thick,-] (0,-1) to (.6,-1) ;
	\draw[thick,-] (.6,-1) to (1,-.6) ;
	\draw[thick,-] (1,-.6) to (1,0) ;
	
	\draw[thick,-] (0,-1) to [out=315,in=225] (.6,-1) ;
	\draw[thick,-] (.6,-1) to [out=0,in=270] (1,-.6) ;
	\draw[thick,-] (1,-.6) to [out=45,in=315] (1,0) ;
	
	\draw[fill] (0,-1) circle [radius=0.035] node[below=2pt]  {$(2,0)$};
	\draw[fill] (.6,-1) circle [radius=0.035] node[right=6pt]  {$(2,1)$};
	\draw[fill] (1,-.6) circle [radius=0.035] node[right=4pt]  {$(2,2)$};
	\draw[fill] (1,0) circle [radius=0.035] node[right=4pt]  {$(2,3)$};
	
	%4th sheet
	\draw[thick,-] (1,1) to (1,1.6) ;
	\draw[thick,-] (1,1.6) to (.6,2) ;
	\draw[thick,-] (.6,2)  to (0,2) ;
	
	\draw[thick,-] (1,1) to [out=45,in=315] (1,1.6) ;
	\draw[thick,-] (1,1.6) to [out=90,in=0] (.6,2) ;
	\draw[thick,-] (.6,2) to [out=135,in=45] (0,2) ;
	
	\draw[fill] (1,1) circle [radius=0.035] node[below right=2pt]  {$(3,0)$};
	\draw[fill] (1,1.6) circle [radius=0.035] node[right=6pt]  {$(3,1)$};
	\draw[fill] (.6,2) circle [radius=0.035] node[right=6pt]  {$(3,2)$};
	\draw[fill] (0,2) circle [radius=0.035] node[above=4pt]  {$(3,3)$};
	
	%Inter-sheet B-edges
	\draw[thick, dotted,-] (-2,1) to (-2,0) ;
	\draw[thick, dotted,-] (-1,-1)  to (0,-1) ;
	\draw[thick, dotted,-] (1,0) to (1,1) ;
	\draw[dashed,-] (0,2) to (-1,2) ;
	%Inter-sheet A-edges
	
	\draw[dashed,-] (-2,1) to (1,1) ;
	\draw[thick, dotted,-] (1,0) to (-2,0) ;
	\draw[thick, dotted,-] (-1,-1) to (-1,2) ;
	\draw[thick, dotted,-] (0,2) to (0,-1) ;
	\draw[-] (-.5,-2) node[]{$\G_2^{BSV} \simeq ~\G_1^{BSV} \g ~\G_1^{BSV}$}; 
	\end{tikzpicture}
	\caption{The Schreier graph $\G_1^{BSV},$ the Tile graph $\G_1'^{BSV}$ of BSV torsion-free group and the graph $\G_1^{BSV} \g \G_1^{BSV}$.}
	\label{Figure:torsion2.4}
\end{figure}
\section{Galois coverings of Schreier graphs}\label{sec:galois}

Our aim in this section is to present the main theorem, which generalizes Proposition~2.1.4 of \cite{shaikh2016zeta}. Recall that $G$ is a group generated by a bounded automaton $\A$, and that the action of $G$ on each level of the rooted tree $X^*$ is transitive, where $X$ is a finite alphabet of cardinality $d$. Let $S$ denote the generating set of $G$, and let $S' \subseteq S$ be the set of reachable states. By Assumption~\ref{ass:bdd}, the set $\pc$ of post-critical sequences is non-empty, and therefore $S'$ is also non-empty.

\begin{mypro}\label{Pro:auto_exist}
Let $\G_{n+1}$ be a $d$-sheeted unramified covering of $\G_n$, and let $e \in E = E(\G_n) \setminus E(\G_n')$. Then there exists an automorphism $\sigma_e : \G_{n+1} \rightarrow \G_{n+1}$ such that 
\begin{equation}
    \sigma_e(wx_i) = w\psi_{s|_u}(x_i), \label{eqn:auto_non-tile}
\end{equation}
for all $w \in X^n$.
\end{mypro}
\begin{proof}
Let $e = \{u, s(u)\} \in E$. Then $u$ is a post-critical vertex in $\G_n$, and $s|_u$ is a reachable state. By Proposition~13.3 of \cite{terras2010zeta}, $e$ admits a unique lift $\tilde{e}$ in $\G_{n+1}$. According to Equation~\eqref{Equation:gr-rot3}, the lift $\tilde{e} = \{ux_i, s(u)s|_u(x_i)\}$ starts on the $x_i$-th sheet and ends at $s(u)s|_u(x_i)$, which lies on the $x_j$-th sheet, where $x_j = s|_u(x_i)$.

Since $s|_u(x_i) = \psi_{s|_u}(x_i)$, where $\psi_{s|_u} \in Sym(X)$ is the root permutation associated with the reachable state $s|_u$, we define a map $\sigma_e : \G_{n+1} \to \G_{n+1}$ by
$$\sigma_e(wx_i) = w\psi_{s|_u}(x_i), \quad \text{for all } w \in X^n.$$
This map $\sigma_e$ permutes the $d$ sheets according to the root permutation $\psi_{s|_u}$, and preserves adjacency by construction, hence it is an automorphism. This completes the proof.
\end{proof}

We now state our main theorem.

\begin{mythm}\label{Theorem:Galois}
Let $G$ be a group generated by a bounded automaton $\A$. Then the following statements are equivalent:
\begin{enumerate}
    \item \label{item:root} 
    There exists a cyclic subgroup $H\le Sym(X)$ of order $d$ such that $\psi_s\in H$ for all reachable states $s\in S'$, and each $\psi_s$ is either $\idperm$ or has order $d$.
    \item \label{item:frob} The graph $\G_{n+1}$ is a $d$-fold cyclic covering of $\G_n$.
\end{enumerate}
\end{mythm}

\begin{proof} Assume \ref{item:root}. By Proposition~\ref{Proposition:coveringmap}, the Schreier graph $\G_{n+1}$ is a $d$--sheeted unramified covering of $\G_n$, with sheets indexed by the alphabet $X=\{x_1,\dots,x_d\}$.

Let $e=\{u,s(u)\}\in E$ be a non-Tile edge. Then $u$ is a post-critical vertex and the restriction $s|_u$ is a reachable state.
By Proposition~\ref{Pro:auto_exist}, there exists an automorphism $\sigma_e\in Aut(\G_{n+1})$ such that $$\sigma_e(w x_i)= w\,\psi_{s|_u}(x_i)\qquad\text{for all } w\in X^n,\ x_i\in X.$$
That is, $\sigma_e$ acts on the set of sheets via the root permutation $\psi_{s|_u}\in Sym(X)$ and satisfies $\pi\circ\sigma_e=\pi$, where
$\pi:\G_{n+1}\to\G_n$ is the covering map.

Since $d>1$ and $\G_{n+1}$ is connected, there exists at least one reachable state $s|_u\in S'$ such that $\psi_{s|_u}$ is nontrivial. By assumption (\ref{item:root}), all root permutations of reachable states lie in a fixed cyclic subgroup $H\le Sym(X)$ of order $d$, and $\psi_{s|_u}$ has order $d$. Hence $\psi_{s|_u}$ is a $d$-cycle and therefore generates $H$.

Consequently, the powers
$$\sigma_e^k(w x_i)= w\,\psi_{s|_u}^k(x_i),\qquad k=1,\dots,d,$$ define a cyclic subgroup of $Aut(\G_{n+1})$ of order $d$, acting
transitively on the set of sheets and satisfying $\pi\circ\sigma_e^k=\pi$ for all $k$. Thus this subgroup is the Galois group of the covering. Therefore, $\G_{n+1}$ is a $d$-fold cyclic covering of $\G_n$.

We shall now prove the converse. Assume that $\G_{n+1}$ is a $d$-fold cyclic covering of $\G_n$, with Galois group $\Gal=Gal(\G_{n+1}|\G_n)$. Let $e=\{u,s(u)\}\in E$ be a non-Tile edge of $\G_n$. Then $u$ is a post-critical vertex and $s|_u$ is a reachable state. By Proposition~\ref{Pro:auto_exist}, the edge $e$ determines an automorphism $\sigma_e\in\Gal$ whose action on the sheets is given by the root permutation $\psi_{s|_u}\in Sym(X)$. 

Since $d>1$ and $\G_{n+1}$ is connected, there exists at least one non-Tile edge $e$ such that $\sigma_e$ is nontrivial. Because $\Gal$ is cyclic of order $d$ and acts transitively on the set of sheets, the action of any generator of $\Gal$ on the sheets is a $d$-cycle. In particular, the action of $\sigma_e$ on the sheets is a $d$-cycle.

The action of $\Gal$ on the sheets induces an injective homomorphism $$\Gal \hookrightarrow Sym(X),$$  whose image is therefore a cyclic subgroup $H\le Sym(X)$ of order $d$. For any reachable state $s\in S'$, the associated root permutation $\psi_s$ coincides with the action on sheets of some $\sigma_e\in\Gal$, and hence $\psi_s\in H$.

It follows that for each reachable state $s\in S'$, the associated root permutation $\psi_s$ lies in the cyclic subgroup $H\le Sym(X)$; moreover, if $\psi_s$ is nontrivial, then it acts transitively on the $d$ sheets and hence has order $d$. Therefore, each $\psi_s$ is either the identity or a $d$--cycle, which establishes (\ref{item:root}).

\end{proof}

\subsection{Examples}
This section is dedicated to examples of both normal (Galois) and non-normal coverings of Schreier graphs of groups generated 
by bounded automata.

\begin{myexa} 
Listed below are some examples of automata groups for which the Schreier graph $\G_{n+1}$ is a $d$-sheeted normal covering of $\G_n$:
\end{myexa}
\begin{enumerate}
	\item \label{Example:Grigorchuk-normal} For the Grigorchuk group (illustrated in Figure~\ref{Figure:auto_grig}), the covering degree is $d = 2$.  
	\item \label{Example:BSV-normal} For the BSV torsion-free $d$-group, the covering degree is $d = 4$.
	\item \label{Example:Basilica-normal} For the Basilica group, the covering degree is $d = 2$.
\end{enumerate} 

\noindent\textbf{Note:} Although not discussed in detail in this paper, for the Gupta--Sidki $p$-group, the covering degree is $d = p$, and for the Fabrykowski--Gupta group, it is $d = 3$.

\begin{myexa} \label{Example:all-nonnormal} 
Non-normal coverings of Schreier graphs for the Tower of Hanoi group.
\end{myexa} Let $X = \{0,1,2\}$. The automorphisms of the ternary rooted tree $X^*$ are defined by the following wreath recursions:
\begin{equation} 
\label{Equation:wr-tower-group}
a = \psi_a(a, \id, \id),\quad b = \psi_b(\id, b, \id),\quad c = \psi_c(\id, \id, c),
\end{equation}
where $\psi_a = (1~2)$, $\psi_b = (0~2)$, and $\psi_c = (0~1)$ in $\mathrm{Sym}(X)$. For further details on the Schreier graphs of this group, see \cite{DAngeli2012counting}.

The Tower of Hanoi group $\mathbb{T}$ is generated by a bounded automaton $\A$ with $|X| = d = 3$. Although all nontrivial states are reachable, the associated root permutations are not $3$-cycles. Consequently, $\mathbb{T}$ does not satisfy condition~\eqref{item:root} of Theorem~\ref{Theorem:Galois}. Hence, the covering $\G_{n+1} \to \G_n$ is a $3$-sheeted non-normal covering for the Tower of Hanoi group.

\section{$L$ functions and zeta functions of graphs}\label{sec:zeta}

The Galois coverings $\G_{n+1} \to \G_n$ arising from groups generated by bounded automata admit a natural decomposition of their cycle structure governed by the action of the Galois group. Ihara zeta functions and their Artin $L$-function refinements provide a convenient framework for organizing this information representation-theoretically. In this section, we exploit the Galois nature of these coverings to express the Ihara zeta function of the covering graph $\Y \to Y$ in terms of Artin $L$-functions associated with irreducible representations of the Galois group $\Gal=Gal(\Y|Y)$. This approach allows explicit computations in our setting and clarifies how automaton-theoretic symmetries are reflected in zeta-function invariants of the Schreier graphs.

Let $\Y = \G_{n+1}$ be a normal covering of $Y = \G_n$. The Ihara zeta function of the covering graph $\Y$ can be expressed in terms of the Artin $L$-functions. To compute these $L$-functions, we require the irreducible representations of the Galois group $\Gal$. Associated with each irreducible representation $\rho$ of $\Gal$, there is an Artin $L$-function. The Ihara zeta function $\zeta_{\Y}(t)$ is given as the product of these Artin $L$-functions.

\medskip
\noindent
\textbf{Convention} for vertices and sheets of a normal cover: Suppose $\Y \to Y$ is a normal covering with Galois group $\Gal$. We designate one sheet of $\Y$ as the identity sheet, denoted sheet $\gid$. For any $\sigma \in \Gal$, the image of the identity sheet under $\sigma$ is referred to as sheet $\sigma$. Every vertex $\tilde{v}$ in $\Y$ can be uniquely written as $\tilde{v} = (v, \sigma)$, where $v = \pi(\tilde{v})$ and $\sigma$ indicates the sheet containing $\tilde{v}$.

\medskip
\begin{mydef} The Artin $L$-function associated with a representation $\rho$ of $\Gal$ is defined by a product over prime cycles in $Y$:
$$
L(u,\rho,\Y|Y) = \prod_{[C]~\text{prime in } Y} \det(I - \rho(\sigma(\widetilde{C})) \, u^{\nu(C)})^{-1},
$$
where $\widetilde{C}$ is the lift of $C$ to the graph $\Y$ starting on sheet $\gid$ and ending on sheet $\tau$. The element $\sigma(\widetilde{C})$ is the Frobenius automorphism defined in Definition 16.3 of \cite{terras2010zeta} as
$$
\sigma(\widetilde{C}) = \tau^{-1},
$$
with $\gid, \tau \in \Gal$. The matrix $\rho(\sigma(\widetilde{C}))$ is a $d_{\rho} \times d_{\rho}$ matrix associated with the irreducible representation $\rho$, where $d_{\rho}$ is the degree of $\rho$.\end{mydef}

\medskip
\begin{mydef} The $m \times m$ matrix $A(\sigma)$ for $\sigma \in \Gal$ is defined such that the $(i,j)$-entry is
\begin{equation}
A(\sigma)_{i,j} = \text{ number of edges in } \Y \text{ between } (i, \gid) \text{ and } (j, \sigma),
\label{Equation:sheet1toandother}
\end{equation}
where $\gid$ denotes the identity element in $\Gal$ and $m$ is the number of vertices in $Y$. \end{mydef}

\medskip
We now recall Proposition 2.1 of \cite{dedeo2005zeta}.

\begin{mypro} \label{block-diagonalization} Let $\Y \to Y$ be a normal covering with Galois group $\Gal$. Then the adjacency matrix of $\Y$ can be block-diagonalized, with each block of the form
\begin{equation}
A_{\rho} = \sum_{\sigma \in \Gal} A(\sigma) \otimes \rho(\sigma),
\label{Equation:blockdiagonal}
\end{equation}
and each block appears $d_{\rho}$ times, where $\rho$ is an irreducible representation of $\Gal$, and $d_{\rho}$ is its degree.\end{mypro}
\medskip
Define $Q_{\rho} = Q \otimes I_{d_{\rho}}$, where $Q$ is the diagonal degree matrix and $I_{d_{\rho}}$ is the identity matrix of size $d_{\rho}$. Then the following analogue of the Ihara determinant formula holds (see Theorem 18.15 of \cite{terras2010zeta}):
$$
L(t,\rho,\Y|Y)^{-1} = (1 - t^2)^{(r - 1)d_{\rho}} \det(I - t A_{\rho} + t^2 Q_{\rho}).
$$

Hence, the zeta function of the covering graph $\Y$ can be expressed as:
\begin{equation}
\label{Equation:zetaasL}
\zeta_{\Y}(t) = \prod_{\rho \in \widehat{\Gal}}  L(t, \rho, \Y|Y)^{d_{\rho}}.
\end{equation}

where the product is over all irreducible representations of $\Gal$.

\medskip
For further details, see Corollary 18.11 on page 155 of \cite{terras2010zeta}. In the next subsection, we will apply this formulation to compute zeta and $L$-functions for examples arising from automaton groups.

\subsection{Examples}

In this section, we explore the Ihara zeta and $L$-functions associated with the Schreier graphs of automaton groups. Through these examples, we aim to illustrate the application of the theoretical concepts discussed previously and to demonstrate how these functions can reveal important structural information about the underlying groups.

\begin{myexa} \label{Example:zeta-BSV} BSV-group \end{myexa}
By Theorem \ref{Theorem:Galois}, the covering $\Y = \G_2^{BSV} \simeq \G_1^{BSV} \g \G_1^{BSV}$ over the graph $Y = \G_1^{BSV}$, as depicted in Figure \ref{Figure:torsion2.4}, represents a $4$-sheeted normal covering. 

To obtain a spanning subgraph of $Y = \G_1^{BSV}$, we remove the edges $e_a = \{3, a(3) = 0\}$ and $e_b = \{3, b(3) = 0\}$. This gives the Tile graph $\G_1'^{BSV}$ of $\G_1^{BSV}$. We visualize the covering graph $\Y$ by placing $4$ sheets $\G_1'^{BSV}$ of $Y$ and labeling each sheet as indicated in Table 2. The connections between the four sheets in the covering graph $\Y$ are outlined in Table 3.

In this scenario, the irreducible characters of the cyclic Galois group $$\Gal = \langle \sigma = (0,1,2,3) \mid \sigma^4 = \gid \rangle \simeq \mathbb{Z}/4\mathbb{Z}$$ are given by $$\chi_i(j) = \exp\left(\frac{2\pi i j}{4}\right)$$ for $1 \leq i, j \leq 4.$

Next, we present the matrices $A(\sigma), \sigma \in \Gal$, as defined in Equation \eqref{Equation:sheet1toandother}:

$$
A(\gid) = \begin{pmatrix}
0 & 2 & 0 & 0 \\
2 & 0 & 2 & 0 \\
0 & 2 & 0 & 2 \\
0 & 0 & 2 & 0 \\
\end{pmatrix}, \quad
A(\sigma) = A(\sigma^3) = \begin{pmatrix}
0 & 0 & 0 & 1 \\
0 & 0 & 0 & 0 \\
0 & 0 & 0 & 0 \\
1 & 0 & 0 & 0 \\
\end{pmatrix}, \quad
A(\sigma^2) = O_{4 \times 4}.
$$

Using Equation \eqref{Equation:blockdiagonal}, we compute the Artinized adjacency matrices $A_{\chi_i}$:

$$
A_{\chi_1} = \begin{pmatrix}
0 & 2 & 0 & 2 \\
2 & 0 & 2 & 0 \\
0 & 2 & 0 & 2 \\
2 & 0 & 2 & 0 \\
\end{pmatrix}, \quad
A_{\chi_2} = A_{\chi_4} = A(\gid), \quad
A_{\chi_3} = \begin{pmatrix}
0 & 2 & 0 & -2 \\
2 & 0 & 2 & 0 \\
0 & 2 & 0 & 2 \\
-2 & 0 & 2 & 0 \\
\end{pmatrix}.
$$

Note that $A_{\chi_1} = A_Y$.

\begin{table}[ht]
	\caption{Notation for sheets.}
	\centering
	\begin{tabular}{c|c|c}
		Vertex set of $\Y$ & Vertex set of $\G_2^{BSV}$ & Group element \\
		\hline \hline
		$\{0,1,2,3\} \sim \{(0,0), (0,1), (0,2), (0,3)\}$ & $\{00, 10, 20, 30\}$ & $\gid$ \\
		$\{0',1',2',3'\} \sim \{(1,0), (1,1), (1,2), (1,3)\}$ & $\{01, 11, 21, 31\}$ & $\sigma = (0123)$ \\
		$\{0^{''},1^{''},2^{''},3^{''}\} \sim \{(2,0),(2,1), (2,2), (2,3)\}$ & $\{02, 12, 22, 32\}$ & $\sigma^2 = (02)(13)$ \\
		$\{0^{'''},1^{'''},2^{'''},3^{'''}\} \sim \{(3,0), (3,1), (3,2), (3,3)\}$ & $\{03, 13, 23, 33\}$ & $\sigma^3 = (0321)$ \\
	\end{tabular}
\end{table}

\begin{table}[ht]
	\caption{Connections between sheet $\gid$ in $\Y$ and other sheets.}
	\centering
	\begin{tabular}{c|c}
		Vertex & Adjacent vertices in $\Y$ \\
		\hline \hline
		$0$ & $1, 1, 3', 3^{'''}$ \\
		$1$ & $0, 0, 2, 2$ \\
		$2$ & $1, 1, 3, 3$ \\
		$3$ & $2, 2, 0', 0^{'''}$ \\
	\end{tabular}
\end{table}

We now determine the reciprocals of the $L$-functions for the covering $\Y|Y$:

1. For $A_{\chi_1}$:
$$
\zeta_{Y}(t)^{-1} = L(t, A_{\chi_1}, \Y|Y)^{-1} = (1 - t^2)^4 (t - 1)(t + 1)(3t - 1)(3t + 1)(3t^2 + 1)^2.
$$

2. As $A_{\chi_2} = A_{\chi_4}$:
$$
L(t, A_{\chi_2}, \Y|Y)^{-1} = L(t, A_{\chi_4}, \Y|Y)^{-1} = (1 - t^2)^4 (9t^4 - 6t^3 + 2t^2 - 2t + 1)(9t^4 + 6t^3 + 2t^2 + 2t + 1).
$$

3. For $A_{\chi_3}$:
$$
L(t, A_{\chi_3}, \Y|Y)^{-1} = (1 - t^2)^4 (9t^4 - 2t^2 + 1)^2.
$$

By Equation \eqref{Equation:zetaasL}, we obtain:
\begin{eqnarray}
\zeta_{\Y}(t)^{-1} &=& (1 - t^2)^{16} (t - 1)(t + 1)(3t - 1)(3t + 1)(3t^2 + 1)^2 (9t^4 - 2t^2 + 1)^2 \nonumber \\
&& \times (9t^4 - 6t^3 + 2t^2 - 2t + 1)^2 (9t^4 + 6t^3 + 2t^2 + 2t + 1)^2. \nonumber
\end{eqnarray}

Note: Equation \eqref{Equation:zetaasL} provides an iterative method for computing $\zeta_{\G_{n+1}}$ from $\zeta_{\G_n}$, given the matrices $A_{\rho}$. Using this method, one can directly compute the Artin $L$-functions and Ihara zeta functions for the Schreier graphs associated with the Fabrykowski–Gupta, Basilica, Grigorchuk, and Gupta–Sidki $p$-groups.

\section{Non-bounded Automaton Groups}\label{sec:nonbounded}

In this section, we study the Schreier graphs of a group generated by a non-bounded automaton.

\begin{myexa} \label{Example:lin-group-ex} Linear Automaton Group $\L$. \end{myexa}
Let $X = \{0,1\}$. Define the automorphisms of the binary tree $X^*$ by the following wreath recursions:
\begin{equation}
\label{Equation:wr-lin-group}
a = \psi_a(a, \id),\quad b = \psi_b (b, a),
\end{equation}
where $\psi_a$ is the transposition $(0~1) \in \mathrm{Sym}(X)$ and $\psi_b = \idperm$.

This group $\L$ is generated by $a$ and $b$. It provides the simplest example of a group generated by an automaton with \emph{polynomial growth}, yet without \emph{bounded activity}. For more details on this group, see \cite{bondarenko2012family}. For illustrations of the Moore diagram of $\L$ and the Schreier graphs $\G_1^{\L}$, $\G_2^{\L}$, and $\G_3^{\L}$ of $\L$, refer to Figures \ref{Figure:lin-group-automaton}, \ref{Figure:L-Sch12}, and \ref{Figure:L-Sch3}, respectively.
\begin{figure}[!htb]
	\begin{tikzpicture}[scale=1.5]
	\tikzstyle{every state}=[draw=black, thick,fill=none]
	\node[state] at (-2,0) (a) {$b$};
	\node[state] at (0,0) (b) {$a$};
	\node[state] at (2,0) (e) {$\id$};
	\draw[thick,->] (a) to (b);
    \node[above] at (-1,0) {$1|1$};
	\path[->] (a) edge [loop above] node[above] {$0|0$} ();
	\draw[thick,->] (b) to node[above] {$1|0$} (e);
	\path[->] (b) edge [loop above] node[above] {$0|1$} ();
	\path[->] (e) edge [loop above] node[above] {$0|0,1|1$} ();
	\end{tikzpicture}
	\centering
	\caption{The Moore diagram of the linear automaton of $\L$}
	\label{Figure:lin-group-automaton}
\end{figure}
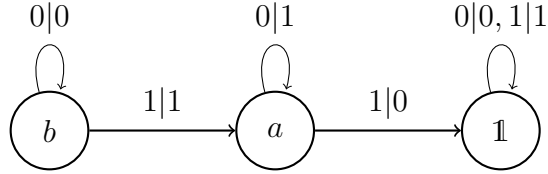
\begin{figure}[!htb]
	\centering
    \begin{subfigure}[b]{0.45\textwidth}
        \begin{tikzpicture}[scale=2.5]
		%($\G_1$)
		\draw[thin] (0,1) to [out=45,in=135] (1,1);
		\draw[thin] (1,1) to [out=225,in=315] (0,1);
		\draw[thin] (1.1,1) circle [radius=0.1];
		\draw[thin] (-0.1,1) circle [radius=0.1];
		\draw[fill] (0,1) circle [radius=0.035]; 
		\draw (0,1) node[below=2pt]  {$0$};
		\draw[fill] (1,1) circle [radius=0.035]; 
		\draw (1,1) node[below=2pt]  {$1$}; 
		\end{tikzpicture}\label{Figure:L-Sch1:left}
        \caption{The graph $\G_1^{\L}$}
    \end{subfigure}
    \begin{subfigure}[b]{0.45\textwidth}
        \begin{tikzpicture}[scale=2]
		\draw[thin] (0,1) to (1,0) to (0,-1) to (-1,0) to (0,1);
		\draw[thin] (-1,0) to  [out=-90,in=180] (0,-1.3) to [out=0,in=-90] (1,0);
		\draw[thin] (-1,0) to  [out=90,in=180] (0,1.3) to [out=0,in=90] (1,0);
		\draw[thin] (0,1.1) circle [radius=0.1];
		\draw[thin] (0,-1.1) circle [radius=0.1];
		\draw[fill] (0,1) circle [radius=0.035]; 
		\draw[fill] (1,0) circle [radius=0.035]; 
		\draw[fill] (0,-1) circle [radius=0.035]; 
		\draw[fill] (-1,0) circle [radius=0.035]; 
		\draw (0,1) node[right=2pt]  {$00$};
		\draw (1,0) node[right=2pt]  {$11$};
		\draw (0,-1) node[right=2pt]  {$01$};
		\draw (-1,0) node[left=2pt]  {$10$};
		\end{tikzpicture}\label{Figure:L-Sch2:right}
    \caption{The graph $\G_2^{\L}$}
    \end{subfigure}
	\caption{The Schreier graphs $\G_1^{\L},\G_2^{\L}$ of $\L$.}
	\label{Figure:L-Sch12}
\end{figure}
\begin{figure}[!htb]
	\centering
	\begin{tikzpicture}[scale=1.3]
	\draw[thin] (0,2) to (1.2,1.2) to (2,0) to (1.2,-1.2) to (0,-2) to (-1.2,-1.2) to (-2,0) to (-1.2,1.2) to (0,2);
	\draw[thin] (0,2.1) circle [radius=0.1];
	\draw[thin] (0,-2.1) circle [radius=0.1];
	\draw[thin] (1.2,1.2) to (1.2,-1.2) to (-1.2,-1.2) to (-1.2,1.2) to (1.2,1.2);
	\draw[thin] (2,0) to [out=90,in=0] (0,2.3) to [out=180,in=90] (-2,0);
	\draw[thin] (2,0) to [out=-90,in=0] (0,-2.3) to [out=180,in=-90] (-2,0);
	\draw[fill] (0,2) circle [radius=0.035]; 
	\draw[fill] (1.2,1.2) circle [radius=0.035]; 
	\draw[fill] (2,0) circle [radius=0.035]; 
	\draw[fill] (1.2,-1.2) circle [radius=0.035]; 
	\draw[fill] (0,-2) circle [radius=0.035]; 
	\draw[fill] (-1.2,-1.2) circle [radius=0.035]; 
	\draw[fill] (-2,0) circle [radius=0.035]; 
	\draw[fill] (-1.2,1.2) circle [radius=0.035]; 
	\draw (0,2) node[right=2pt]  {$000$};
	\draw (1.2,1.2) node[below left]  {$111$};
	\draw (2,0) node[right=2pt]  {$011$};
	\draw (1.2,-1.2) node[above left]  {$101$};
	\draw (0,-2) node[right=2pt]  {$001$};
	\draw (-1.2,-1.2) node[above right]  {$110$};
	\draw (-2,0) node[left=2pt]  {$010$};
	\draw (-1.2,1.2) node[below right]  {$100$};
	\end{tikzpicture}
	\caption{The Schreier graph $\G_3^{\L}$ of $\L$.}
	\label{Figure:L-Sch3}
\end{figure}
\begin{myexa}\label{linear-non-normal} Galois Coverings of Schreier Graphs of $\L$. \end{myexa}
To construct the spanning tree $T_n$ of $\G_n^{\L}$, let $V(T_n) = X^n$ and define $u_0 = \{0\}^n$ as the word of length $n$ consisting entirely of zeros. For an element $g \in\, \L$, $g^{2^k}$ denotes $g$ composed with itself $2^k$ times. An edge $e = \{u,v\}$ exists in $\G_n^{\L}$ if there is a generator $s \in S$ such that $s(u) = v$. We refer to such an edge as the $s$-edge at $u$. The action of $a$ on $X^n$ has order $2^n$, which generates the spanning tree $T_n$ of $\G_n^{\L}$:

\begin{eqnarray}
\label{edges}
E(T_n) &=& \bigl\{ \{a^k(u_0),\, a^{k+1}(u_0)\} \;\big|\; 0 \le k \le 2^n-2 \bigr\}, \\
\label{vertices}
V(T_n) &=& \bigl\{ a^k(u_0) \;\big|\; 0 \le k \le 2^n-1 \bigr\}.
\end{eqnarray}

\begin{mypro}\label{pro:linear}
Let $\G_m^{\L}$ and $\G_n^{\L}$ be the level-$m$ and level-$n$ Schreier graphs of $G$. Then the following hold:
\begin{enumerate}
  \item \label{item:unrami_L}
  If $m>n$, then $\G_m^{\L}$ is an unramified covering of $\G_n^{\L}$.
  \item \label{item:galois_L}
  If $m=n+1$, then $\G_m^{\L}$ is a $2$-fold cyclic covering of $\G_n^{\L}$.
\end{enumerate}
\end{mypro}

\begin{proof}
\begin{enumerate}
  \item
  Let $m=n+k$ with $k\ge1$. Define a map
  $$\pi:\G_m^{\L}\longrightarrow \G_n^{\L}, \qquad
  \pi(x_1\cdots x_m)=x_1\cdots x_n,$$
  which forgets the last $k$ letters.
  The map $\pi$ is the composition of the successive projections
 $$ \G_{n+k}^{\L}\to \G_{n+k-1}^{\L}\to \cdots \to \G_n^{\L},$$
  each of which is an unramified covering. Hence $\pi$ is itself an unramified covering.
  The covering is unramified since local degrees are preserved under each projection.

  \item
  Let $T_n$ be a spanning tree of $\G_n^{\L}$. Using Equation~\eqref{vertices}, define
  
  $$T_n0=\{a^i(u_0)0 \mid 0\le i\le 2^n-1\}, \qquad
  T_n1=\{a^i(u_0)1 \mid 0\le i\le 2^n-1\}.$$
  
  These sets form a partition of the preimage of $T_n$ in $\G_{n+1}^{\L}$ and hence
  represent the two sheets of the covering.

  Define a map $\sigma:\G_{n+1}^{\L}\to \G_{n+1}^{\L}$ by
  $$ \sigma(x_1\cdots x_n y)=a^{2^n}(x_1\cdots x_n y).$$
  Since $a^{2^n}$ fixes the prefix $x_1\cdots x_n$ and acts on the last letter via the
  root permutation $\psi_a$, we have
 $$ \sigma(x_1\cdots x_n y)=x_1\cdots x_n\,\psi_a(y).$$
  As $\psi_a$ is a transposition, $\sigma$ has order $2$. Moreover,
  $\pi\circ\sigma=\pi$. Therefore, $\G_{n+1}^{\L}\to\G_n^{\L}$ is a $2$-fold cyclic covering.
\end{enumerate}
\end{proof}

\subsection{Open Question}

This study motivates the following intriguing problem. Let $G$ be a group generated by a \emph{non-bounded} automaton $\A$ acting transitively on a finite alphabet $X$ with $|X| = d$. Let $S' \subset S$ be the set of reachable states, and let $\psi_s$ denote the root permutation of each reachable $s \in S'$. This leads to the following question:

\medskip
\noindent
{\bf Question:}  
Let $G$ be a group generated by a non-bounded automaton $\A$. Suppose that for each state $s \in S'$, the root permutation $\psi_s$ is either the identity or a $d$-cycle in $\mathrm{Sym}(X)$.  
Does this condition imply that, $\G_{n+1}$ is a $d$-fold cyclic covering of $\G_n$?  
If not, what additional hypotheses (if any) recover an analogue of the bounded case?

\end{document}